\documentclass[twocolumn]{autart}    

\usepackage{graphicx}          
\usepackage{amsmath}
\usepackage{amssymb}
\usepackage{amsfonts}
\usepackage{color}
\usepackage[T1]{fontenc}
\newtheorem{lemma}{Lemma}
\newtheorem{Remark}{Remark}
\newtheorem{proposition}{Proposition}
\newtheorem{corollary}{Corollary}
\newtheorem{definition}{Definition}

\def\qed{ \rule{.08in}{.08in}}
\newcommand{\F}{\mathbb{F}}
\newcommand{\N}{\mathcal{N}_{\rm out}}
\newcommand{\Nin}{\mathcal{N}_{\rm in}}
\newcommand{\I}{\mathcal{I}}

\begin{document}

\begin{frontmatter}

\title{Stability Structures of Conjunctive Boolean Networks\thanksref{footnoteinfo}} 

\thanks[footnoteinfo]{This paper was not presented at any IFAC meeting. This research was supported in part by AFOSR MURI Grant FA 9550-10-1-0573.  }

\author[UIUC]{Zuguang Gao}\ead{zgao19@illinois.edu},    
\author[CUB]{Xudong Chen}\ead{xudong.chen@colorado.edu},               
\author[UIUC]{Tamer Ba\c{s}ar}\ead{basar1@illinois.edu}  

\address[UIUC]{University of Illinois at Urbana-Champaign, United States}  
\address[CUB]{University of Colorado at Boulder, United States}             


\begin{abstract}                   
A Boolean network is a finite dynamical system, whose variables take values from a binary set. {\color{black}The value update rule for each variable is a Boolean function, depending on a selected subset of variables.} Boolean networks have been widely used in modeling {\color{black}gene regulatory networks}. We focus in this paper on a special class of Boolean networks, termed as \emph{conjunctive Boolean networks}. A Boolean network is \emph{conjunctive} if the associated value update rule is comprised of only AND operations. It is known that any trajectory of a finite dynamical system will enter a periodic orbit. We characterize in this paper all periodic orbits of a conjunctive Boolean network {\color{black}whose underlying graph is strongly connected}. In particular, we establish a bijection between the set of periodic orbits and the set of binary necklaces of a certain length. We further investigate the stability of a periodic orbit. Specifically, we perturb a state in the periodic orbit by changing the value of a single entry of the state. The trajectory, with the perturbed state being the initial condition, will enter another (possibly the same) periodic orbit in finite time steps. We then provide a complete characterization of all such transitions from one periodic orbit to another. In particular, we construct a digraph, with the vertices being the periodic orbits, and the (directed) edges representing the transitions among the orbits. We call such a digraph the stability structure of the conjunctive Boolean network. 
\end{abstract}

\end{frontmatter}

\section{Introduction}
Finite dynamical systems are discrete-time dynamical systems with finite state spaces. They have a long and successful history of being used in biological networks~\cite{funahashi1993approximation}, epidemic networks~\cite{khanafer2014stability}, social networks~\cite{etesami2016complexity}, and engineering control systems~\cite{imer2006optimal}. 
In this paper, we focus on a special class of finite dynamical systems, called Boolean networks (or Boolean automata networks~\cite{noual2013non}).  Boolean networks are finite dynamical systems whose variables are of Boolean type, usually labeled as ``$1$'' and ``$0$''. {\color{black}The Boolean function, also known as the value update rule, for each variable depends on a selected subset of the variables. }   

Boolean networks have been widely used  in systems biology and (mathematical) computational biology. 
This line of research began with Boolean network representations of molecular networks~\cite{kauffman1969homeostasis}, and was later generalized to the so-called logical models~\cite{thomas1990biological}.
Since then there have been some studies of various classes of Boolean functions which are particularly suited to the logical expression of gene regulation~\cite{kauffman1969metabolic,thomas1973boolean,raeymaekers2002dynamics}. Evidence has been provided in~\cite{sontag2008effect} that biochemical networks are ``close to monotone''. Roughly speaking, a Boolean network is monotonic if its Boolean function has the property that the output value of the function for each variable is non-decreasing if the number of ``$1$''s in the inputs increases.  
Monotonic Boolean networks have been studied both theoretically~\cite{jarrah2010dynamics,melliti2013convergence,noual2012updating,noual2012boolean,remy2003description} and in applications~\cite{georgescu2008gene,mendoza1999genetic}. 
Also, there have been studies of Boolean networks with other types of Boolean functions: For example, Boolean networks whose Boolean functions are monomials were studied in~\cite{colon2005boolean,colon2006monomial,park2014monomial}. The work by~\cite{veliz2010dynamics} considers the dynamics of the systems where the Boolean functions are comprised of semilattice operators, i.e., operators that are commutative, associative, and idempotent. Boolean networks whose Boolean functions are comprised of only XOR operations were investigated in~\cite{alcolei2015flora}, and whose Boolean functions are comprised of AND and NOT operations were studied in~\cite{veliz2012and,veliz2015dimension}. 

A special class of Boolean functions, of particular interest to us here, is the so-called nested canalyzing functions. This class of functions was introduced in~\cite{stuart1993origins}, and often used to model genetic networks \cite{harris2002model,kauffman2003random,kauffman2004genetic}. Roughly speaking, a canalyzing function is such that if an input of the function holds a certain value, called the ``canalyzing value'', then the output value of the function is uniquely determined regardless of the other values of the inputs~\cite{jarrah2007nested}. The majority of Boolean functions that appear in the literature on Boolean networks are nested canalyzing functions. Among the nested canalyzing functions, there are two simple but important classes: A function in the first class is comprised of only AND operations, with ``$0$'' the canalyzing value,  while a function in the second class is comprised of only OR operations, with ``$1$'' the canalyzing value. 
The corresponding Boolean networks are said to be {\em conjunctive}  and {\em disjunctive}, respectively~\cite{jarrah2010dynamics,goles2012disjunctive}. Note that there is a natural isomorphism between the class of conjunctive Boolean networks and the class of disjunctive Boolean network: indeed, if $f$ (resp. $g$) is a function on $n$ Boolean variables $x_1,\ldots, x_n$, comprised of only AND (resp. OR) operations, then 
$
f(x_1,\ldots, x_n) = \neg g(\neg x_1,\ldots, \neg x_n),
$
where ``$\neg$'' is the negation  operator, i.e., $\neg 0 = 1$ and $\neg 1 = 0$.  
It thus suffices to consider only conjunctive Boolean networks. We note here that a conjunctive Boolean network is monotonic.

Since a Boolean network is a finite dynamical system, for any initial condition, the trajectory generated by the system will enter a periodic orbit (also known as a limit cycle) in finite time steps {\color{black}(see, for example,~\cite{colon2005boolean})}. A question that comes up naturally is how the dynamical system behaves if a ``perturbation'' occurs in a state of a periodic orbit---meaning that one (and only one) of the variables fails to follow the update rule for the next time step (a precise definition is given in Subsection~\ref{stabilitystruc}). The trajectory, with the perturbed state as its initial condition, will then enter another periodic orbit (possibly return to the original orbit). One of the questions addressed in this work is thus to characterize all possible transitions among the periodic orbits upon the occurrence of a perturbation.

A complete characterization of these transitions among the periodic orbits is given in Theorem~\ref{thm1}, which captures the stability structure of a conjunctive Boolean network. 
The analysis of Theorem~\ref{thm1} relies on a representation of periodic orbits, which identifies the orbits with the so-called binary necklaces (a definition is given in Subsection~\ref{binarynecklace}). In particular, we show that there is a bijection between the set of periodic orbits and the set of binary necklaces of a certain length.  To establish this bijection, we introduce in Section~III a new approach for analyzing the system behavior of a conjunctive Boolean network: Roughly speaking, we decompose the original Boolean network into several components. For each of the components, there corresponds an induced dynamics. We then relate in Theorem~\ref{dynamic} the original dynamic to these induced dynamics and establish several necessary and sufficient conditions for a state to be in a periodic orbit. This new approach may be of independent interest as it can be applied to other types of Boolean networks as well.      

The rest of the paper is organized as follows.
In Section~\ref{pre}, we first provide some basic definitions and notations  for directed graphs and the binary necklace. We then introduce the class of conjunctive Boolean networks in precise terms. Some preliminary results on such networks are also given.  
In Section~\ref{decomposing}, we introduce the new approach as mentioned above. A detailed organization will be given at the beginning of that section. Then, in Section~\ref{stability}, we characterize all possible transitions among periodic orbits. Moreover, we associate each transition with a positive real number, termed as transition weight, which can be understood as the likelihood of the occurrence of the transition.   
We provide conclusions and outlooks in Section~\ref{end}. The paper ends with an Appendix which contains proofs of some technical results.

\section{Preliminaries}\label{pre}

\subsection{Directed graph}
We introduce here some useful notation associated with a directed graph (or simply digraph). Let $D=(V,E)$ be a directed graph. We denote by $v_iv_j$ an edge from $v_i$ to $v_j$ in $D$. We say that $v_i$ is an {\em in-neighbor} of  $v_j$ and $v_j$ is an {\em out-neighbor} of $v_i$. 
The sets of in-neighbors and out-neighbors of vertex $v_i$ are denoted by $\mathcal{N}_{\rm in}(v_i)$ and
$\mathcal{N}_{{\rm out}}(v_i)$, respectively.
The  \emph{in-degree} and \emph{out-degree}
of vertex $v_i$ are defined as $|\mathcal{N}_{{\rm in}}(v_i)|$ and $|\mathcal{N}_{{\rm out}}(v_i)|$, respectively.  

Let $v_i$ and $v_j$ be two vertices of $D$. A {\em walk} from $v_i$ to $v_j$, denoted by $w_{ij}$, is a sequence $v_{i_0}v_{i_2}\cdots v_{i_m}$ (with $v_{i_0} = v_i$ and $v_{i_m} = v_j$) in which $v_{i_k}v_{i_{k+1}}$ is an edge of $D$ for all 
$k\in\{0,1,\ldots,m-1\}$. A walk is said to be a {\it path} if all the vertices in the walk are pairwise distinct. A {\em closed walk} is a walk $w_{ij}$ such that the starting vertex and ending vertex are the same, i.e., $v_i=v_j$.
A walk is said to be a {\em cycle} if there is no repetition of vertices in the walk other than the repetition of the starting- and ending-vertex. The {\it length} of a path/cycle/walk is defined to be the number of edges in that path/cycle/walk. 


%
A \emph{strongly connected graph} is a directed graph such that for any two vertices $v_i$ and $v_j$ in the graph, there is a path from $v_i$ to $v_j$. A \emph{cycle digraph} is a directed graph that consists of a single cycle.

\subsection{Binary necklace}\label{binarynecklace}
\begin{figure}[h]
	\centering
	\includegraphics[height=40mm]{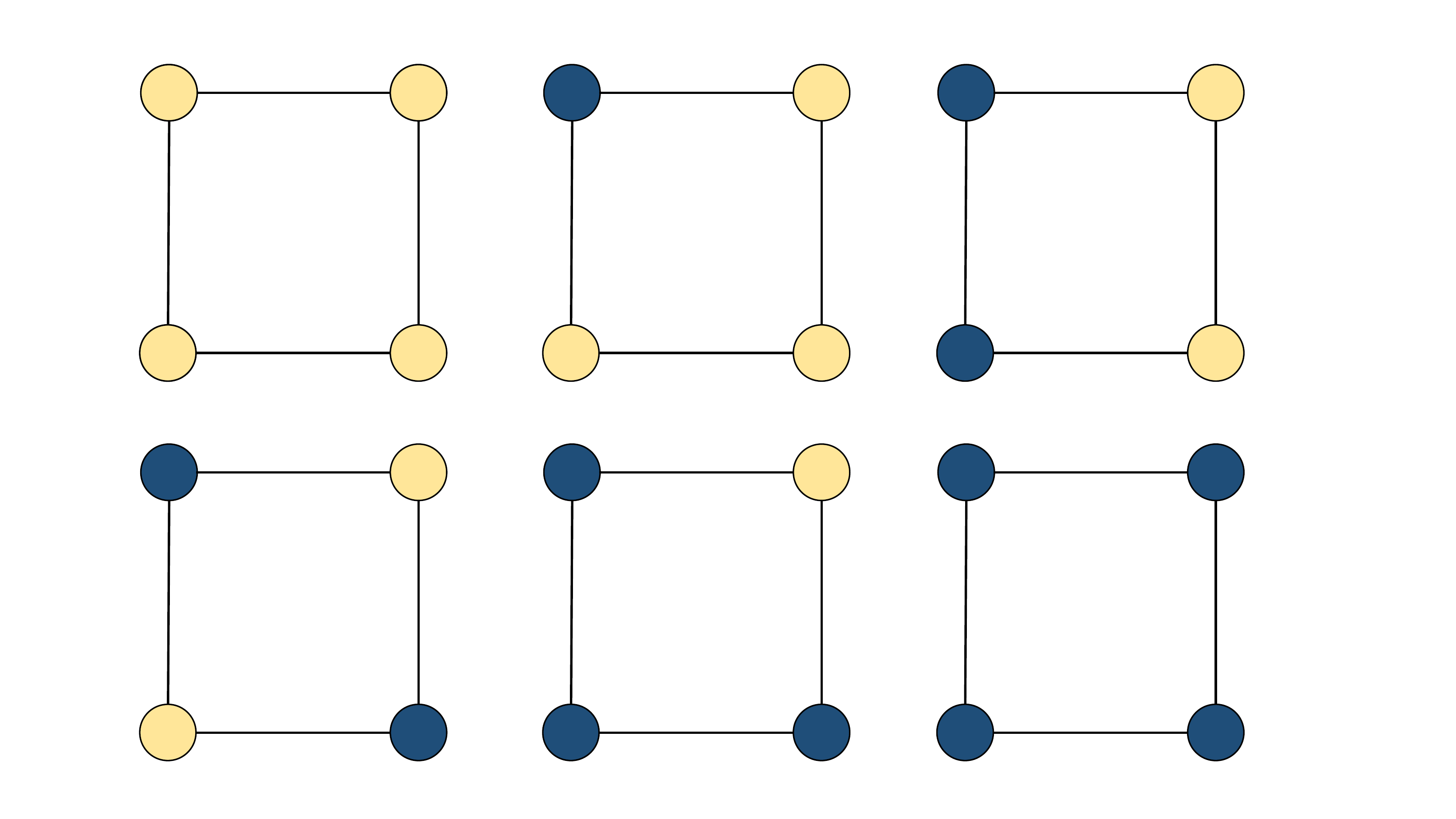}
	\caption{All binary necklaces of length~$4$. If the bead is plotted in dark blue (resp. light yellow), then it holds value ``1'' (resp, ``0''). }
	\label{necklace}
\end{figure}

A {\bf binary necklace} of length $p$ is an equivalence class of $p$-character strings over the binary set $\mathbb{F}_2=\{0,1\}$, taking all rotations circular shifts as equivalent. For example, in the case of {\color{black}$p = 4$}, there are six different binary necklaces, as illustrated in Fig.~\ref{necklace}. 
A \emph{necklace with fixed density} is a necklace in which the number of zeros (and hence, ones) is fixed. 
The {\bf order} of a necklace is the cardinality of the corresponding equivalence class, and it is always a divisor of $p$.  
An \emph{aperiodic necklace} (see, for example,~\cite{varadarajan1990aperiodic}) is a necklace of order~$p$, i.e., no two distinct rotations of a necklace from such a class are equal. Thus, an aperiodic necklace cannot be partitioned into more than one sub-strings which have the same alphabet pattern. For example, a necklace of $1010$ (row 2, column 1 in Fig.~\ref{necklace}) can be partitioned into two substrings $10$ and $10$ which have the same alphabet pattern, and thus is not aperiodic. A necklace of $1000$ (row 1, column 2 in Fig.~\ref{necklace}) cannot be partitioned into more than one sub-strings with the same alphabet pattern, and is aperiodic.

\subsection{Conjunctive Boolean network}

Let $\mathbb{F}_2=\{0,1\}$ be the finite field with two elements. The two elements ``$0$'' and $``1''$ can, for example, represent the ``off'' status and ``on'' status of a gene, respectively. 
We call a function $g$ on $n$ variables a {\em Boolean function} if it is of the form $g:\F^n_2 \to \F_2$.  The so-called \emph{Boolean network} on $n$ Boolean variables $x_1(t),\ldots, x_n(t)$ is a discrete-time dynamical system, whose update rule can be described by a set of Boolean functions $f_1,\ldots, f_n$: 
$$
x_i(t+1) = f_i(x_1(t),\ldots, x_n(t)), \hspace{5pt} \forall i = 1, \ldots, n.  
$$
For convenience, we let $x(t):= (x_1(t),\ldots, x_n(t)) \in \F_2^n$ be the state of the Boolean network at time~$t$. We further let $$f:= (f_1, \ldots, f_n): x(t) \mapsto x(t+1).$$ We refer to $f$ as the {\bf value update rule} associated with the Boolean network. Note that following this value update rule, all Boolean variables {\color{black}update} their values synchronously (in parallel) at each time step. We refer to~\cite{goles2012disjunctive,noual2013non,ruz2013preservation} for results on Boolean networks with asynchronous (sequential) updating schemes. 

It is well known that for any initial condition $x(0)\in \F_2$, the trajectory $x(0), x(1), \ldots$ will enter a {\em periodic orbit} in a finite amount time. Precisely, there exists a time $t_0 \ge 0$ and an integer number $p \ge 1$ such that $x(t_0 + p) = x(t_0)$. Moreover, if $x(t_0 +q) \neq x(t_0)$ for any $q = 1,\ldots, p-1$, then the sequence $\{x(t_0), \ldots, x(t_0 + p-1)\}$, taking rotations as equivalent, is said to be a {\bf periodic orbit}, and we call $p$ its {\bf period}. If the period of a periodic orbit is one, then $x(t_0) = x(t_0 + 1) = \cdots$. We then call the state $x(t_0)$ a {\bf fixed point}.


We consider, in this paper,  a special class of Boolean networks, termed {\em conjunctive Boolean networks}. Roughly speaking, a Boolean network is conjunctive if each Boolean function $f_i$ is an AND operation on a selected subset of the $n$ variables. We provide below a precise definition:

\begin{definition}[Conjunctive Boolean network~\cite{jarrah2010dynamics}]\label{CBN}
	A Boolean network $f = (f_1, \ldots, f_n)$  is {\bf conjunctive} if each Boolean function $f_i$, for all $i = 1,\ldots, n$, can be expressed as follows:
	\begin{equation}\label{eq:updaterule}
	f_i(x_1,\ldots, x_n)=\prod^n_{j = 1} x_j^{\epsilon_{ji}}
	\end{equation}
	with $\epsilon_{ji}\in \{0,1\}$ for all $j=1,\ldots,n$. 
\end{definition}

Note that if we let $I_i:= \{j \mid \epsilon_{ji} = 1\}$, then $f_i$ is nothing but an AND operator on the variables $x_j$, for $j\in I_i$. 


We now associate a conjunctive Boolean network with a \emph{dependency graph}, whose definition is given as follows.


\begin{definition}[Dependency graph~\cite{jarrah2010dynamics}]
	Let $f = (f_1, \ldots, f_n)$ be the value update rule associated with a conjunctive Boolean network. The associated {\bf dependency graph} is a directed graph $D = (V,E)$ of $n$ vertices. An edge from $v_i$ to $v_j$, denoted by $v_iv_j$,  exists in $E$ if and only if $i\in I_j$.   
\end{definition}

\begin{Remark}
	A conjunctive Boolean network uniquely determines its dependency graph. Conversely, given a digraph $D$, there is a unique conjunctive Boolean network whose dependency graph is $D$.  
\end{Remark}


In this paper we assume that the dependency graph $D$ is strongly connected. 
We now present some preliminary results on the network and the associated digraph. 

First, note that if a digraph $D = (V, E)$ is strongly connected, then it can be written as the union of its cycles (\cite{chen2015consensus}): Let $D_1 = (V_1, E_1),\ldots, D_N = (V_N, E_N)$, with $V_i \subset V$ and $E_i\subset E$, be the cycles of $D$. Then, 
$$
D = \left (\cup^N_{i = 1} V_i, \cup^N_{i = 1} E_i \right ).
$$
Said in another way, each vertex of $D$ is contained in at least one cycle of $D$. Now, let $n_i$ be the length of $D_i$. Then, we have the following fact for the possible periods of the conjunctive Boolean network:

\begin{lemma}
	A positive integer $p$ is the period of a periodic orbit of a conjunctive Boolean network if and only if $p$ divides the length of each cycle.
	\label{relation39}
\end{lemma}

\begin{Remark}
	Note that if the greatest common divisor of the cycle lengths is one, then the period $p$ of a periodic orbit $\{x(t_0), \ldots, x(t_0+p-1)\}$ has to be one, and hence $x(t_0)$ is a fixed point of the conjunctive Boolean network. 
\end{Remark}

We refer to~\cite{jarrah2010dynamics,gao2016periodic} for proofs of Lemma~\ref{relation39}. We further have the following fact:

\begin{lemma}\label{lem:fixedpoint}
	A state $x\in \F^n_2$ is a fixed point of a conjunctive Boolean network if and only if all the $x_i$'s hold the same value.
\end{lemma}

\noindent
{\bf Proof.}
	It should be clear that if all the $x_i$'s hold the same value, then $x$ is a fixed point. We now show that the converse is also true. The proof is done by contradiction: assume that there are two vertices $v_i$ and $v_j$ such that $x_i = 0$ and $x_j =1$. Since the dependency graph $D$ is strongly connected, there is a walk $w_{ij}$ from $v_i$ to $v_j$. Let $l(w_{ij}) = q$, and label the vertices along the walk as follows:
	$
	w_{ij} = v_{k_0}v_{k_1}\ldots v_{k_q},
	$
	with $v_{k_0} = v_{i}$ and $v_{k_q} = v_j$. Now, suppose that $x_{k_0}(t_0) = 0$; then, from~\eqref{eq:updaterule}, we have $x_{k_1}(t_0+1) = 0$, and $x_{k_2}(t_0+2) = 0,\ldots, x_{k_q}(t_0 + q) = 0$. On the other hand, since $x$ is a fixed point, $x_{k_q}(t_0 + q) = x_{k_q}(t_0) = 1$, which is a contradiction.  
	\hfill$\qed$


\section{Irreducible Components of Strongly Connected Graphs} \label{decomposing}
Let $D = (V,E)$ be the dependency graph associated with a conjunctive Boolean network. Assume that $D$ is strongly connected, and recall that $D_1,\ldots, D_N$ are cycles of $D$, and $n_1,\ldots, n_N$ are their lengths. Now, let $p^*$ be the greatest common divisor of $n_i$, for $i = 1,\ldots, N$: $$p^*:= {\rm gcd}\{n_1, n_2, ..., n_N\}.$$ This is also known as the \emph{loop number} of $D$~\cite{colon2005boolean}.  
The digraph $D$ is said to be {\bf irreducible} if $p^*=1$. If the digraph $D$ is not irreducible, then we show in this section that there is a decomposition of $D$ into $p^*$ components each of which is irreducible. This section is thus organized as follows: In Subsection~\ref{VSP}, we partition the vertex set $V$ in a particular way into $p^*$ subsets. Following this partition, we then construct, in Subsection~\ref{irrcomp}, $p^*$ digraphs, as we call the irreducible components of $D$, whose vertex sets are the $p^*$ partitioned subsets. We show in Proposition~\ref{stronglyirreducible} that each irreducible component is indeed irreducible, and moreover, strongly connected.  Then, in Subsection~\ref{def:indyn}, we define a conjunctive Boolean network, as we call an induced dynamics, on each irreducible component. We further establish the relationships between the original dynamics and the $p^*$ induced dynamics.

\subsection{Vertex set partition}\label{VSP}


Following Lemma~\ref{relation39}, we introduce a partition of the vertex set $V$.  Roughly speaking, the partition is defined such that  the vertices in a partitioned subset are connected by walks whose lengths are multiples of a common divisor of the cycle lengths.  
We now define the partition in precise terms. To proceed, we first have some definitions and notations.  
Let $v_i,v_j$ be any two vertices in $V$, and $w_{ij}$ be a walk from  $v_i$ to $v_j$. We denote by $l(w_{ij})$ the length of $w_{ij}$. 

\begin{definition}\label{rel}
	Let $p$ divide the lengths of cycles of the dependency graph~$D$. We say that a vertex $v_i$ is {\bf related to} another vertex $v_j$ (or simply write $v_i\sim_p v_j$) if there exists a walk $w_{ij}$ from $v_i$ to $v_j$ such that $p$ divides $l(w_{ij})$.
\end{definition}

We note here that the relation introduced in Definition~\ref{rel} is in fact an equivalence relation. Specifically, we have the following fact:  

\begin{lemma}\label{equivrel}
	The relation $\sim_p$ is an equivalence relation, i.e., for any $v_i,v_j,v_r\in V$, the following three properties hold:
	\begin{enumerate}
		\item Reflexivity: $v_i\sim_p v_i$. 
		\item Symmetry: $v_j\sim_p v_i$ if and only if $v_i\sim_p v_j$. 
		\item Transitivity: if $v_i\sim_p v_j$ and $v_j\sim_p v_r$, then $v_i\sim_p v_r$. 
	\end{enumerate}\,
\end{lemma}

We refer to the Appendix for a proof of the lemma. With the preliminaries above, we construct  a subset of $V$ as follows: First, choose an arbitrary vertex $v_i$ as a base vertex; then, we define
\begin{equation}\label{eq:defWpi}
[v_i]_p:=\{v_j\in V\mid v_j\sim_p v_i\}.
\end{equation}
Note that from Lemma~\ref{equivrel}, the subset $[v_i]_p$, for any $v_i\in V$, is an equivalence class of $v_i$. 
We further establish the following result: 
\begin{proposition}\label{417}
	The following two properties hold: 
	\begin{enumerate}
		\item If $v_i\sim_p v_j$, then $[v_i]_p=[v_j]_p$. If $v_i\not\sim_p v_j$, then $[v_i]_p\cap[v_j]_p=\varnothing$.
		\item Let $v_0\in V$, and choose vertices $v_1,\ldots, v_{p-1}$ such that $v_1\in \mathcal{N}_{\rm out}(v_0), \ldots, v_{p-1}\in \mathcal{N}_{\rm out}(v_{p-2})$. Then, the subsets $[v_0]_p,\ldots, [v_{p-1}]_p$ form a partition of $V$:
		\begin{equation}
		V = \sqcup^{p-1}_{i = 0}  [v_i]_p.
		\label{eq:partition}
		\end{equation}
	\end{enumerate}\,
\end{proposition}

\begin{figure}[h]
	\centering
	\includegraphics[height=60mm]{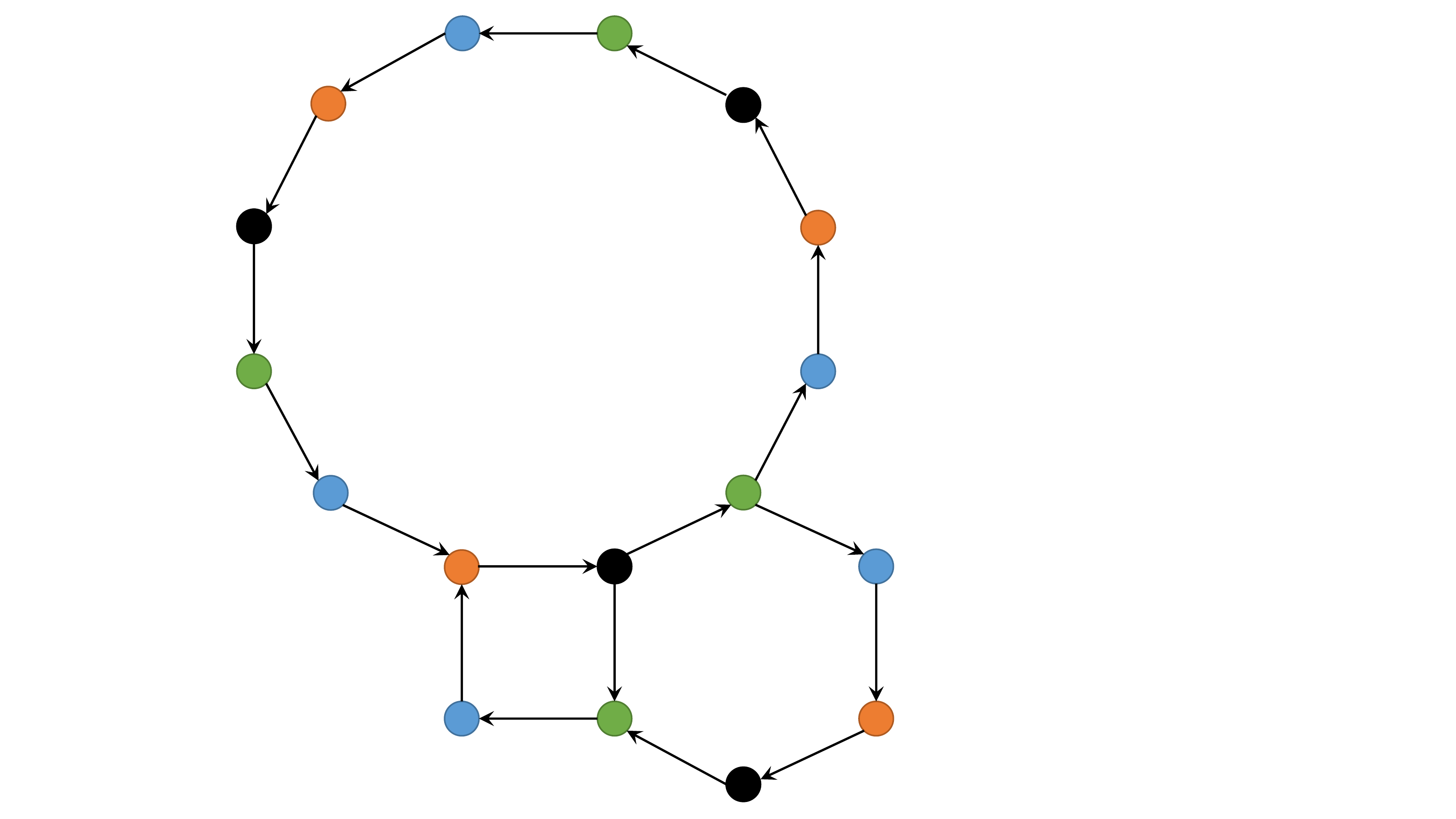}
	\caption{The digraph in the figure has three cycles whose lengths are $4$, $8$, and $12$, respectively. Let $p = 4$ be a common divisor of the cycle lengths. Then, the associated partition yields $4$ disjoint subsets, with the vertices of the same color belonging to the same subset.  
	}
	\label{P}
\end{figure}

We provide in Fig.~\ref{P} an example of such a partition of~$V$.  

The remainder of {\color{black}this} subsection is devoted to the proof of Proposition~\ref{417}. The first item of the proposition directly follows from the fact that each $[v_i]_p$, for any $v_i \in V$, is an equivalence class of $v_i$. We now prove the second item of the proposition. To proceed, first note that in~\eqref{eq:defWpi}, if $v_j\in [v_i]_p$, then there is a walk $w_{ji}$ from $v_j$ to $v_i$ with $l(w_{ji})$ a multiple of $p$. We now show that if $w_{ji}$ is a walk from $v_j$ to $v_i$, then $l(w_{ji})$ has to be a multiple of $p$.

\begin{lemma}
	Let $v_i\sim_p v_j$, and $w_{ij}$ be an arbitrary walk from $v_i$ to $v_j$. Then, $l(w_{ij})$ is a multiple of $p$.
	\label{relation48}
\end{lemma}

\noindent
{\bf Proof.}
	Since $v_i\sim_p v_j$, there exists a walk $w_{ij}$ such that $l(w_{ij})=k_1p$ for some $k_1\in\mathbb{Z}^+$. Suppose there is a different walk $w'_{ij}$ which connects $v_i$ to $v_j$. We need to prove that $l(w'_{ij})$ is a multiple of $p$. By Lemma~\ref{equivrel} (reflexivity), we have $v_j\sim_p v_i$. Therefore, there exists a walk $w_{ji}$ whose length $l(w_{ji})=k_2p$. {\color{black}Concatenating} $w'_{ij}$ and $w_{ji}$, we get a closed walk $w'_{ii}$.  It is known that in strongly connected graphs, any closed walk can be decomposed into cycles. Since $p$ divides the lengths of all cycles, we have that $p$ divides $l(w'_{ii})$. Now we have that $p$ divides both $l(w_{ji})$ and $l(w'_{ii})$. Thus, $p$ also divides $l(w'_{ii})-l(w_{ji})=l(w'_{ij})$.
	\hfill$\qed$

With Lemma~\ref{relation48} at hand, we are now in a position to complete the proof of Proposition~\ref{417}:

\noindent
{\bf Proof.}
	We prove here item 2 of Proposition~\ref{417}. 
	We first show that the subsets $[v_i]_p$ for $i = 1,\ldots, p$, are pairwise disjoint, and then show that their union is $V$.
	
	Choose a pair $(i,j)$ with $0\leq i<j\leq p-1$. Then, it should be clear that there is a walk $w_{ij}$  from $v_i$ to $v_j$ with $l(w_{ij}) = j - i < p$. So, by Lemma~\ref{relation48}, $v_i\not\sim_p v_j$, and hence $[v_i]_p\cap[v_j]_p=\varnothing$.

	It now suffices to show that $V = \cup^{p-1}_{i = 0}[v_i]_p$. Picking an arbitrary vertex $v_r$, we show that $v_r\in [v_i]_p$ for some $i = 0,\ldots, p-1$. Since the digraph $D$ is strongly connected, there is a walk $w_{r0}$ from $v_r$ to $v_0$. We then write $l(w_{r0}) = kp + q$, with $0\le q \le p-1$. If $q = 0$, then $v_r \in [v_0]_p$. We thus assume that $q \neq 0$. Now, let $w_{0, p-q}$ be a walk from $v_0$ to $v_{p-q}$ with $l(w_{0,p-q}) = p - q$. Then, by concatenating $w_{r0}$ with $w_{0,p-q}$, we obtain a walk $w_{r,p-q}$ from $v_r$ to $v_{p-q}$ with $l(w_{r,p-q}) = (k+1)p$. Thus, $v_r\in [v_{p-q}]_p$. 
	\hfill$\qed$


\subsection{Irreducible components}\label{irrcomp}
Let $D = (V, E)$ be a strongly connected digraph, and $p^*$ be its loop number. {\color{black} For a vertex~$v\in V$, we simply write $[v_0]$ instead of $[v_0]_{p}$ if  $p = p^*$.}
We now decompose the digraph $D$ into $p^*$ components:  

\begin{definition}[Irreducible components]\label{Def:IrrComp}
	Let $D = (V, E)$ be a strongly connected digraph, and $p^*$ be its loop number. {\color{black} Choose a vertex~$v_0$ of $D$, and let $v_1\in \N(v_0),\ldots, v_{p^*-1} \in \N(v_{p^*-2})$.} The subsets  $[v_0],\ldots, [v_{p^*-1}]$ then form a partition of $V$. The {\bf irreducible components} of $D$ are digraphs $G_0 = (U_0,F_0),\ldots, G_{p^*-1} = (U_{p^*-1}, F_{p^*-1})$, with their vertex sets $U_k$'s given by
	$$
	U_k := [v_k], \hspace{10pt} \forall k = 0,\ldots, p^*-1.
	$$
	The edge set $F_k$ of $G_k$ is determined as follows: Let $u_i$ and $u_j$ be two vertices of $G_k$. Then, $u_iu_j$ is an edge of $G_k$ if there is a walk $w_{ij}$ from $u_i$ to $u_j$ in $D$ with $l(w_{ij}) = p^*$. 
\end{definition}

\noindent
{\color{black} As we will see in Proposition~\ref{stronglyirreducible}, each irreducible component of $D$ is indeed irreducible. }

\begin{Remark}
	The walk $w_{ij}$ in the definition above is either a path or a cycle (which is the case if $u_i = u_j$) because otherwise there will be a cycle of $D$ properly contained in $w_{ij}$ which contradicts the fact that $l(w_{ij}) = p^*$ divides all the cycle lengths. If $w_{ij}$ is a cycle, then the edge $u_iu_j$ is a self-loop.    
\end{Remark}

\begin{figure}[h]
	\centering
	\includegraphics[height=60mm]{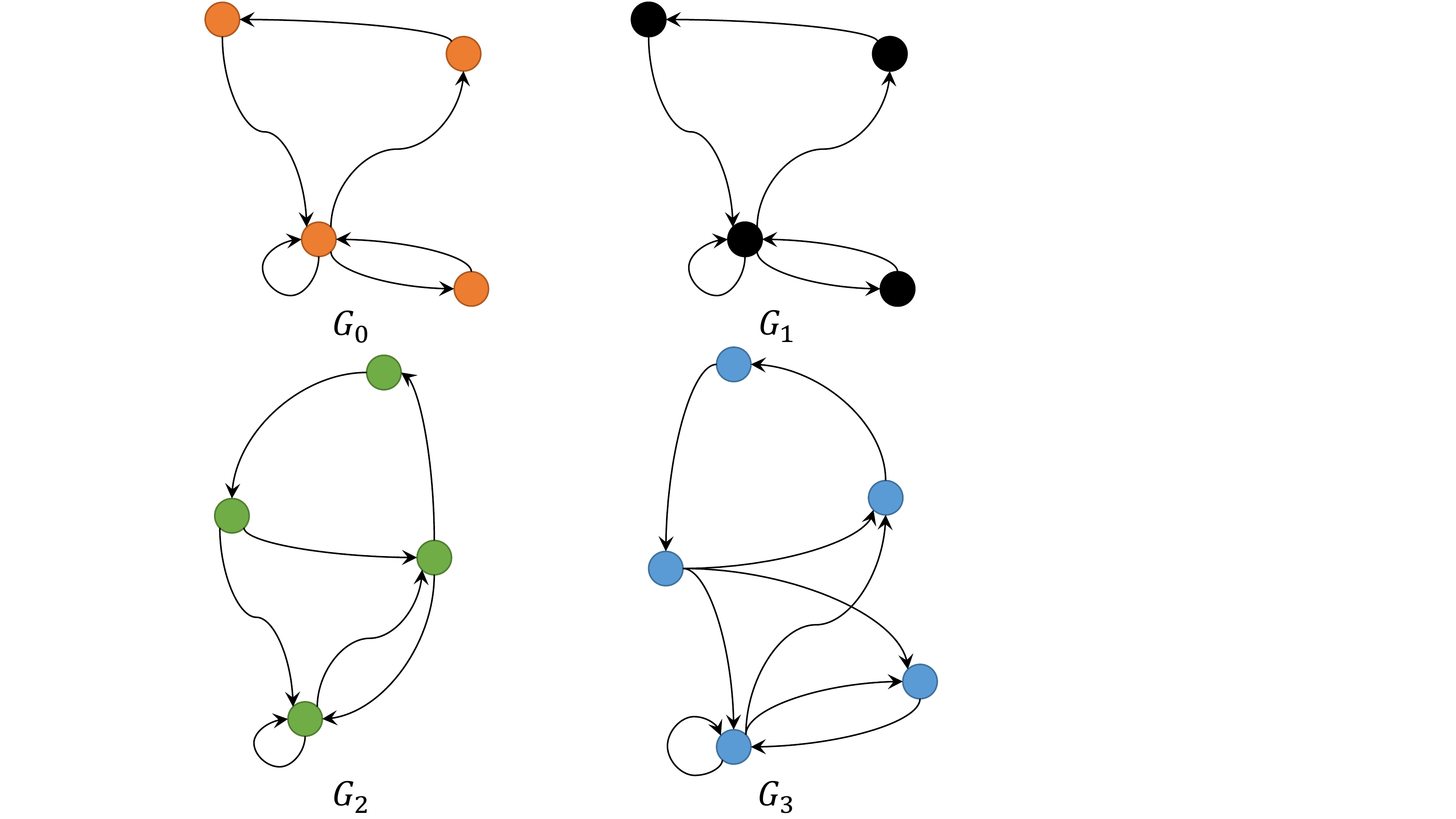}
	\caption{Irreducible components of the digraph shown in Fig.~\ref{P}.}
	\label{components}
\end{figure}

We provide an example in Fig.~\ref{components} in which we show the irreducible components of the digraph shown in Fig.~\ref{P}.

We now establish some properties associated with the irreducible components. We first have the following result:

\begin{proposition}\label{stronglyirreducible}
	Each  $G_k$, for $k=0,\ldots,p^*-1$, is strongly connected and irreducible.
\end{proposition}

To establish the proposition, we need the following lemma: 

\begin{lemma}\label{cycle}
	If there is a cycle of length $n_i$ in the digraph $D$, then there is a cycle of length $n_i/p^*$ in any one of its irreducible components.
\end{lemma}

\noindent
{\bf Proof.}
	Let $D_i = (V_i, E_i)$ be a cycle of length $n_i$ in $D$. 
	Note that $p^*$ divides $|V_i|$. By Definition~\ref{Def:IrrComp}, each irreducible component $G_k = (U_k, F_k)$ contains $n_i/p^*$ vertices of $D_i$. 
	For ease of notation, we let $m:= n_i/ p^*$.  Let 
	$
	V_i \cap U_k = \{u_1,\ldots, u_{m}\}.
	$
	We can further assume that there exist walks $w_{i,i+1}$, for $i = 1,\ldots, m$, from $u_i$ to $u_{i+1}$ in $D$ with $l(w_{i,i+1}) = p^*$ (if $i = m$, we identify $u_{m+1} = u_1$).   
	It then follows from Definition~\ref{Def:IrrComp} that  $u_{1}u_{2}, \ldots, u_{m}u_1$ are edges of $G_k$. Thus, the vertices $u_1,\ldots, u_m$, together with the edges  $u_{1}u_{2}, \ldots, u_{m}u_1$, form a cycle in $G_k$, whose length is ${n_i}/{p^*}$. 
	\hfill$\qed$

	\begin{Remark}\label{R4}
		
		We note here that the converse of Lemma~\ref{cycle} does not hold, i.e., even if there is a cycle of length $m$ in each irreducible component $G_k$, the original digraph $D$ does not necessarily have a cycle of length $mp^*$. A counter example is provided in the Appendix.
		
\end{Remark}

With Lemma~\ref{cycle} at hand, we now prove Proposition~\ref{stronglyirreducible}: 

\noindent
{\bf Proof of Proposition~\ref{stronglyirreducible}.}
	We first prove that each $G_k$ is strongly connected. Let $u_{i}$ and $u_{j}$ be two vertices of $U_k$. We show that there exists a walk in $G_k$ from $u_{i}$ to $u_{j}$. Since $u_i\in [u_j]$, from~\eqref{eq:defWpi}, there is a walk $w_{ij}$ in $D$ with $l(w_{ij}) = r p^*$ for some positive integer~$r$. For a later purpose, we label the vertices, along the walk, as 
	$
	w_{ij} = v_{0}v_1\ldots v_{rp^*}, 
	$ 
	with $v_0 = u_i$ and $v_{rp^*} = u_j$. It then follows from Definition~\ref{Def:IrrComp} that  $v_0,v_{p^*},\ldots, v_{rp^*}$ are vertices of $G_k$. Moreover, $v_0v_{p^*},\ldots, v_{(r-1)p^*}v_{rp^*}$ are edges of $G_k$. So, there is a walk $v_0v_{p^*}\ldots v_{rp^*}$ from $v_0$ to $v_{rp^*}$ in $G_k$.
	
	We next show that $G_k$ is irreducible. 
	Let $D_i=(V_i,E_i)$ be a cycle in $G_k$, with $n_i$ the length of $D_i$. Then, from Lemma~\ref{cycle}, there is a cycle in each $G_k$ whose length is ${n_i}/{p^*}$. 
	We thus conclude that the loop number of each $G_k$ is at most
	$
	{\rm gcd}\{n_1/p^*,\ldots, n_N/p^*\} = 1,
	$
	and hence each $G_k$ is irreducible. 
	\hfill$\qed$

Given a subset $V'$ of $V$ and a nonnegative integer $p$, we define a subset $\Nin^p(V')$ by induction: For $p = 0$, let $\Nin^0(V') := V'$; for $p \ge 1$, we define 
\begin{equation}\label{eq:defninp}
\Nin^p(V') := \cup_{v_j\in \Nin^{p-1}(V')} \Nin(v_j). 
\end{equation}
Similarly, we define $\N^p(V')$ by replacing $\Nin$ with $\N$ in~\eqref{eq:defninp}. With the notations above, we have the following result about the relationships between the vertex sets of the irreducible components:

\begin{proposition}\label{pro:NinNout}
	For $k\ge 0$, we have 
	$$
	\left\{
	\begin{array}{l}
	\N^k(U_0) = U_{(k \bmod p^*)}, \vspace{3pt}\\
	\Nin^k(U_0) = U_{(-k \bmod p^*)}.
	\end{array}
	\right.
	$$\,
\end{proposition}

\noindent
{\bf Proof.}
	We prove here only the first relation $\N^k(U_0) = U_{(k \bmod p^*)}$, the other relation can be established in a similar way. It suffices to show that for any $k = 0,\ldots, p^*-1$, we have $\N(U_k) = U_{(k +1 \bmod p^*)}$. There are two cases: 
	
	{\em Case I: $0\le k \le p^* -2 $}. We first show that $\N(U_k) \subseteq U_{(k +1 \bmod p^*)}$. Let $u \in U_k = [v_k]$, and $u' \in \N(u)$.  Since $D$ is strongly connected, there is a walk $w'$ from $u'$ to $v_k$. Moreover, from Lemma~\ref{relation48}, $l(w') \equiv p^* -1 \bmod p^*$. To see this, note that by concatenating the edge $uu'$ with $w'$, we obtain a walk $w$ from $u$ to $v_k$. Since $u\sim_{p^*} v_k$, $l(w) = l(w') +1 \equiv 0 \bmod p^*.$ Now, using the fact that $v_{k+1}$ is the out-neighbor of $v_k$, we  obtain a walk $w^*$ from $u'$ to $v_{k+1}$ by concatenating $w'$ with the edge $v_{k}v_{k+1}$. Since $l(w^*)$ is a multiple of $p^*$, $u'\in [v_{k+1}] = U_{k+1}$. We now show that $\N(U_k) \supseteq U_{(k +1 \bmod p^*)}$. Let $u' \in U_{k+1}$, and $w'$ be a walk from $u'$ to $v_k$. Then, by the same argument, $l(w') \equiv p^* - 1 \bmod p^*$. Now, let $u\in \Nin(u')$. Then, by concatenating the edge $uu'$ with $w'$, we obtain a walk $w$ from $u$ to $v_k$. Moreover, $l(w)$ is a multiple of $p^*$, and hence $u\in [v_k]$, which implies that $u' \in \N(u) \subseteq \N(U_k)$.

	{\em Case II: $k = p^*-1$}. Let $v_{p^*} \in \N(v_{p^*-1})$. It should be clear that $[v_{p^*}] = [v_0]$. On the other hand, we can apply the arguments above, and obtain that $\N(U_{p^*-1}) = [v_{p^*}]$. So, $\N(U_{p^*-1}) = U_0$. 
	\hfill$\qed$

In the end of this subsection, we introduce a special class of digraphs as follows: 
\begin{definition}\label{rose}
	A digraph $D$ is a {\bf rose} if all the cycles of  $D$ satisfy the following two conditions:
	\begin{enumerate}
		\item They have the same length.
		\item They share at least one common vertex of $D$.
	\end{enumerate}\,
\end{definition}

We provide in Fig.~\ref{rosegraph} an example of a rose. We now have the following result: 

\begin{figure}[h]
	\centering
	\includegraphics[height=30mm]{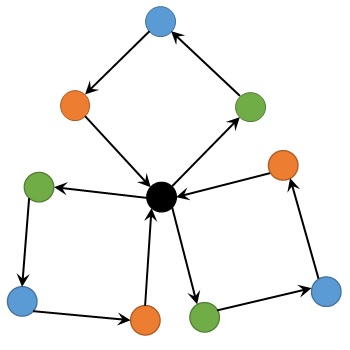}
	\caption{A rose with three cycles of length $4$. The vertex in black is a common  vertex of {\color{black}all} cycles. }
	\label{rosegraph}
\end{figure}


\begin{proposition}\label{special}
	Let $G_k = (U_k,F_k)$, for $k = 0,\ldots,p^*-1$, be irreducible components of $D$. Then, the following hold:
	\begin{enumerate} 
		\item $D$ is a rose if and only if there is at least one $k  \in \{0,\ldots, p^*-1\}$ such that $|U_k| = 1$.
		\item $D$ is a cycle digraph if and only if $|U_k| = 1$, for all $k = 0,\ldots, p^*-1$.
	\end{enumerate}\,
\end{proposition}

We refer to the Appendix for a proof of Proposition~\ref{special}.

\subsection{Induced dynamics}\label{def:indyn}
Let $f = (f_1,\ldots, f_n)$ be a conjunctive Boolean network, and $D$ be the dependency graph. Let $G_0,\ldots, G_{p^*-1}$ be the irreducible components of $D$. Now, for each $k = 0,\ldots, p^*-1$, we can define a conjunctive Boolean network as follows:

\begin{definition}[Induced dynamics]
	An {\bf induced dynamics} on $G_k$ is a conjunctive Boolean network whose dependency graph is $G_k$. 
\end{definition}

We can express the induced dynamics on $G_k$ explicitly as follows: Let $U_k = \{u_1,\ldots, u_m\}$, and $(y_1,\ldots, y_m)$ be the state of the network. Let $g_k= (g_{k_1},\ldots,g_{k_m})$ be the associated value update rule. Then,
$$
g_{k_i}(y_1,\ldots, y_m) = \prod_{u_j\in U_k} y^{\epsilon_{ji}}_j 
$$    
where $\epsilon_{ji} = 1$ if $u_j$ is an in-neighbor of $u_i$ and $\epsilon_{ji} = 0$ otherwise.

We now relate the original dynamics $f$ on $D$ to the induced dynamics on the irreducible components.  We first introduce some notations. Let $V'$ be a subset of $V$. We define $f_{V'}$ to be the restriction of $f$ to $V'$. For a positive integer~$p$, we let $f^p$ be the map defined by applying the map $f$ $p$ times. Given a state $x \in \F^n_2$ and a subset $V'$ of $V$, we let $x_{V'}$ be the restriction of $x$ to $V'$.   
We now establish the main result of this section as follows:

\begin{thm}\label{dynamic}
	Let $G_k = (U_k, F_k)$ be an irreducible component of $D$. Then, the following hold: 
	\begin{enumerate}
		\item Let $g_{k}$ be the induced dynamics on $G_k$. Then, 
		$$
		g_k(x_{U_k}) = f^{p^*}_{U_k}(x), \hspace{10pt} \forall x \in \F^n_2. 
		$$
		\item Suppose that $x(t_0)$ is in a periodic orbit; then,
		\begin{equation}\label{eq:valueshift}
		x_{U_{(k+1 \bmod p^*)}}(t_0 + 1) = x_{U_{k}}(t_0)
		\end{equation}
	\end{enumerate}\,
\end{thm}

We note here that if $x(t_0)$ is in a periodic orbit, then for each $k = 0,\ldots, p^*-1$, the entries of $x_{U_k}(t_0)$ hold the same value. This indeed follows from the first item of Theorem~\ref{dynamic}: 

\begin{corollary}\label{same}
	Let $D = (V, E)$ be the dependency graph of a conjunctive Boolean network, and  $G_k = (U_k, F_k)$, for $k = 0,\ldots,p^*-1$, be its irreducible components. 
	A state $x\in \F^n_2$ is in a periodic orbit of the conjunctive Boolean network if and only if for each $k=0,\ldots,p^*-1$, the entries of $x_{U_k}$ hold the same value.
\end{corollary}

\noindent
{\bf Proof.}
Let $x\in \F^n_2$ be a state. If for each $k = 0,\ldots, p^*-1$, the entries of $x_{U_k}$ hold the same value, then from the first item of Theorem~\ref{dynamic}, $f^{p^*}_{U_k}(x) = g_k(x_{U_k}) = x_{U_k}, $ and hence $f^{p^*}(x) = x$. 
	Conversely, if $x$ is in a periodic orbit of period $p$, then 
	$
	x_{U_k}=f_{U_k}^{p}(x)=f_{U_k}^{p^*}(x)=g_k(x_{U_k}). 
	$
	The first equality holds because $x = f^p(x)$. The second equality holds because $p$ divides $p^*$ (from Lemma~\ref{relation39}). The third equality follows from the first item of Theorem~\ref{dynamic}. 
	So, $x_{U_k}$ is a fixed point of the induced dynamic on $G_k$. From Lemma~\ref{lem:fixedpoint}, we conclude that the entries of $x_{U_k}$ hold the same value.
	\hfill$\qed$

\begin{figure}[h]
	\centering
	\includegraphics[height=40mm]{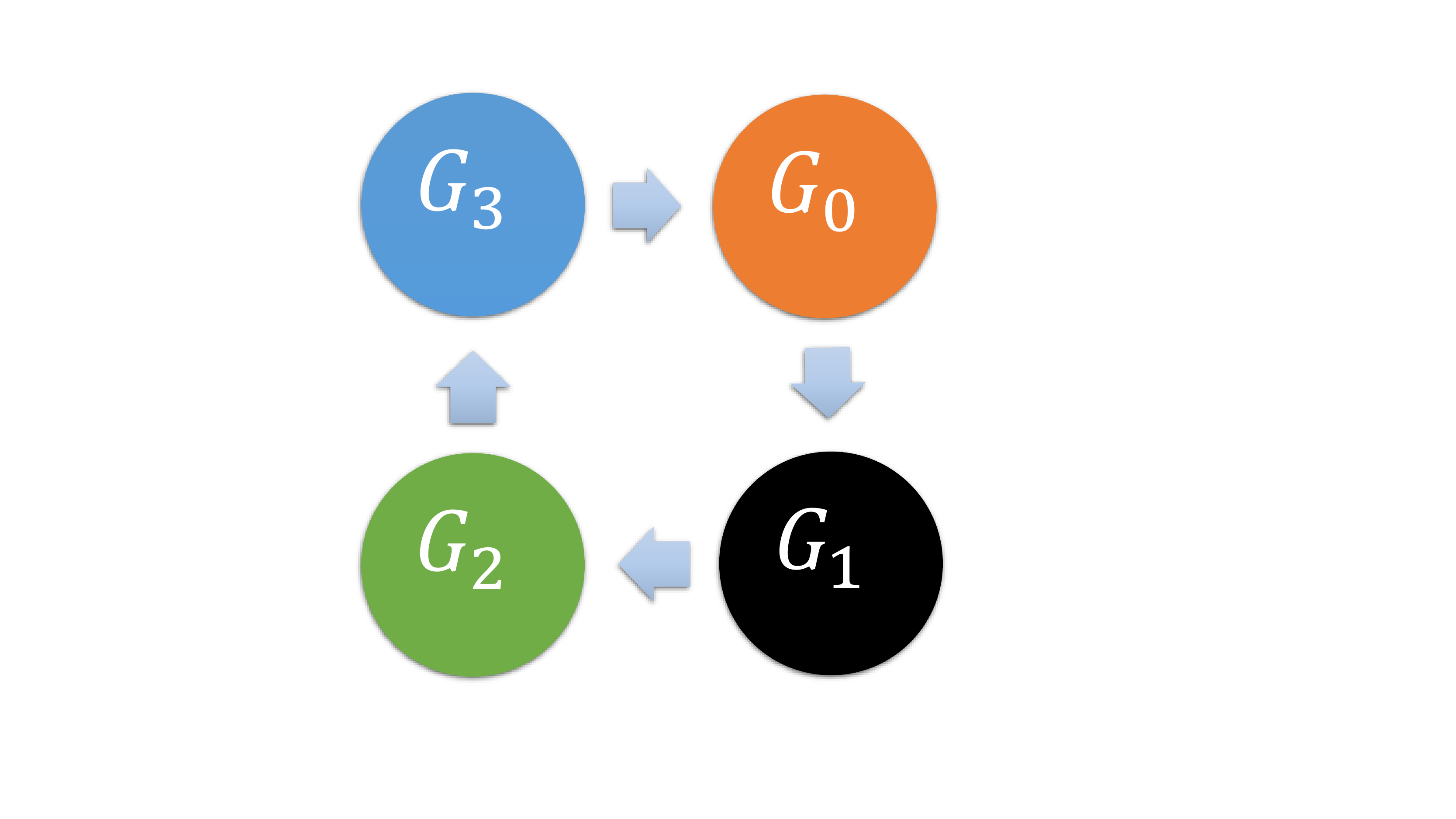}
	\caption{In this figure, $G_0,\ldots, G_3$ are irreducible components of the digraph shown in Fig.~\ref{P}. If $x(t_0)$ is in a periodic orbit, then the vertices of each $G_k$, for $k= 0,\ldots,3$, hold the same value $y_k(t_0)$. Moreover, the value $y_k(t_0)$ will be passed to the vertices of $G_{(k+1) \bmod 4}$ at time step $(t_0 + 1)$.}
	\label{outline}
\end{figure}

So, if $x(t_0)$ is in a periodic orbit, then from the second item of Theorem~\ref{dynamic} and Corollary~\ref{same}, the entries of $x_{U_k}(t_0)$ hold the same value, and moreover, this value will be passed onto the entries of $x_{U_{(k+1 \bmod p^*)}}$ at the next time step. We also illustrate this fact in Fig.~\ref{outline}. 

The remainder of {\color{black}this} section is devoted to the proof of Theorem~\ref{dynamic}.  
For a vertex $v_i$ of $D$ and a positive integer $p$, we define a subset $\Nin^p(v_i)$ of $V$ via induction: For $p = 1$, $\Nin^1(v_i)$ is simply the in-neighbor of $v_i$. For $p \ge 1$, we define $$\Nin^p(v_i):= \cup_{v_j \in \Nin^{p-1}(v_i)} \Nin(v_j).$$ 
In particular, if $p = p^*$ and $v_i$ is a vertex of $G_k$, then from Definition~\ref{Def:IrrComp}, $\Nin^p(v_i)$ is the set of in-neighbors of $v_i$ in $G_k$. 
We further note the following fact:

\begin{lemma}\label{inneigh}
	For any positive integer $p$, we have 
	\begin{equation}\label{eq:poweroff}
	f^p_i(x)= \prod_{v_j \in \Nin^p(v_i)} x_j.
	\end{equation}\,     
\end{lemma}

\noindent
{\bf Proof.}
	We prove the lemma by induction on $p$. For $p=1$, $f_i(x)= \prod_{v_j \in \Nin(v_i)} x_j$, which directly follows from Definition~\ref{CBN}. Now, we assume that~\eqref{eq:poweroff} holds for $p-1$, and prove for $p$. By the induction hypothesis, we have that $x_i(t+p)=\prod_{v_j \in \Nin^{p-1}(v_i)} x_j(t+1).$ From the value update rule, we have that $x_j(t+1)=\prod_{v_k \in \Nin(v_j)} x_k(t).$ So, $x_i(t+p)=\prod_{v_j \in \Nin^{p-1}(v_i)} \prod_{v_k \in \Nin(v_j)} x_k(t).$
	Using the fact that $\Nin^p(v_i) = \cup_{v_j \in \Nin^{p-1}(v_i)} \Nin(v_j),$
	we conclude that~\eqref{eq:poweroff} holds for~$p$.
	\hfill$\qed$

We now prove Theorem~\ref{dynamic}:


\noindent
{\bf Proof of Theorem~\ref{dynamic}.}
	The first item of Theorem~\ref{dynamic} directly follows from Lemma~\ref{inneigh}. 
	We prove here the second item. From the proof of Proposition~\ref{pro:NinNout}, $$U_{(k+1 \bmod p^*)} = \N(U_{k}).$$ 
	So, by the value update rule, $x_{U_{(k+1 \bmod p^*)}}(t_0 + 1)$ depends only on $x_{U_k}(t_0)$.  
	Since $x(t_0)$ is in a periodic orbit, from Corollary~\ref{same}, the entries of $x_{U_{k}}(t_0)$ hold the same value, which then implies that~\eqref{eq:valueshift} holds.
	\hfill$\qed$




\section{Stability of Periodic Orbits}\label{stability}

\subsection{Labeling periodic orbits}\label{Identify}

In this subsection, we find and label all the periodic orbits of a conjunctive Boolean network. Let $D = (V, E)$ be the associated dependency graph, and $p^*$ be its loop number.   Recall that a binary necklace of length $p^*$ is an equivalence class of $p^*$-character strings over $\F_2$, taking rotations as equivalent. The order of a necklace is the cardinality of the equivalence class.

We now show that each periodic orbit can be uniquely identified with a binary necklace of length $p^*$: 
Let $\{x(t_0),\ldots, x(t_0 + p-1)\}$ be a periodic orbit of period $p$. Let $G_k = (U_k, F_k)$, for $k = 0,\ldots, p^*-1$, be the irreducible components of $D$.   
From Corollary~\ref{same}, for each $k = 0,\ldots, p^*-1$, the entries of $x_{U_k}(t_0)$ hold the same value. 
We label these values as $y_0(t_0),\ldots, y_{p^*-1}(t_0)$, with $y_k(t_0)$ being the value of the entries of $x_{U_k}(t_0)$. From the second item of Theorem~\ref{dynamic}, we have that
$$
y_k(t_0 + q) = y_{(k - q \bmod p^*)}(t_0) 
$$ 
for all $k = 0,\ldots, p^* - 1$ and for all $q \ge 0$. This then implies that the periodic orbit $\{x(t_0), \ldots, x(t_0 + p -1)\}$ can be represented by a binary necklace $y_0(t_0)\ldots y_{p^*-1}(t_0)$ whose order is $p$. Conversely, given a binary necklace $y_0\ldots y_{p^*-1}$ of order $p$, we can construct a periodic orbit of period $p$ as follows: Define a state $x \in \F_2^n$ such that the entries of $x_{U_k}$ hold the value $y_k$ for all $k= 0,\ldots, p^*-1$. Appealing again to Corollary~\ref{same} and the second item of Theorem~\ref{dynamic}, we have that $\{x, f(x), \ldots, f^{p-1}(x)\}$ is a periodic orbit of period $p$.  
The arguments above thus imply the following fact: 

\begin{proposition}\label{prop:5}
	There is a bijection between the set of periodic orbits and the set of binary necklaces of length $p^*$. Moreover, such a bijection maps a periodic orbit of period $p$ to a necklace of order $p$.
\end{proposition}

\begin{Remark}\label{R5}
	From the proposition, if two dependency graphs share the same loop number, then the associated conjunctive Boolean networks have the same number of periodic orbits. \end{Remark}

For the remainder of the paper, we let $S$ denote the set of periodic orbits. Each periodic orbit $s\in S$ can be identified with a binary necklace $s = y_0\ldots y_{p^*-1}$. To proceed, we introduce some definitions and notations. Let $\sigma(s)$ be the number of ``$1$''s in the string $s = y_0\ldots y_{p^* - 1}$. We then partition the set $S$ into $(p^* + 1)$ subsets $S_0,\ldots, S_{p^*+1}$:
$$
S_d := \{s\in S\mid \sigma(s) = d\}.
$$ 

Recall that the so-called Euler's totient function $\phi(k)$  counts the total number of integers in the range $[1,k]$ that are relatively prime to $k$. We now present some known results about counting the number of periodic orbits in $S$.

\begin{lemma}\label{old}
	The following two relations hold:
	\begin{enumerate}
		\item For a divisor $p$ of $p^*$, we let $p = \prod^r_{i = 1}p_i^{k_i}$ be its prime factorization. Then, the number of periodic orbits of period $p$ is given by
		$$
		\frac{1}{p}\sum_{i_1=0}^{1}\cdots\sum_{i_r=0}^{1} \left ( (-1)^{\sum^r_{j = 1}i_j} \prod^r_{j = 1} 2^{p_j^{k_j-i_j}} \right ).
		$$
		\item For a number $d = 0,\ldots, p^*$, we have
		$$
		|S_d|=\frac{1}{p^*}\sum_{k\mid {\rm gcd}(p^*-d,d)}\phi(k)\left (\frac{(p^*/k)!}{((p^*-d)/k)!(d/k)!} \right ). 
		$$ 
	\end{enumerate}\,
\end{lemma}

We note here that {\color{black}item~1 of Lemma~\ref{old}} is equivalent to Moreau's necklace-counting formula~\cite{moreau1872permutations}. We refer to~\cite{jarrah2010dynamics} for a proof of item~1, and~\cite{gilbert1961symmetry,ruskey1999efficient} for proofs of item~2.

\subsection{Stability structure}\label{stabilitystruc}

We investigate in this subsection the stability of each periodic orbit of a conjunctive Boolean network. The motivation for this work comes from the fact that the actual process of gene expression is highly complicated. Though conjunctive Boolean networks provide a good model to determine whether a gene can be expressed or not, there are still exceptions and unknown mechanisms that could possibly affect the expression process. We thus want to explore how the system behaves when one gene is not expressed although all necessary proteins are present, or it is expressed even in lack of some necessary proteins.

Let $x(t_0)$ be a  state in a periodic orbit~$s$. We say that a {\em perturbation} occurs at $(t_0 + 1)$ if there is one (and only one) $i\in \{1,\ldots, n\}$ such that $x_{i}(t_0 + 1) = \neg f_i(x(t_0))$. 
As a consequence, $x(t_0 + 1)$ {\color{black} may not be in the periodic orbit~$s$ anymore}. {\color{black} However, }after finite time steps, the system, with $x(t_0+1)$ as its initial condition, will enter a periodic orbit, denoted by~$s'$,  which may or may not be the same as $s$. Our goal in this subsection is to characterize all these {\color{black} transition} pairs $(s,s')$. 

To proceed, we first introduce some definitions and notation. Given a state $x\in \F^n_2$, we let $\mathcal{I}(x)\subset \F^n_2$ be defined as follows: a state $x'$ is in $\mathcal{I}(x)$ if and only if $x'$ differs from $x$ by only one entry, i.e., there is an~$i \in \{1,\ldots, n\}$ such that $x'_i \neq x_i$ and $x'_j = x_j$ for all $j \neq i$. Note that if  $x=x(t_0)$ for $x(t_0)$ a state in a periodic orbit, then $\I(x)$ is the set of states upon the condition that a perturbation occurs at $(t_0 + 1)$. We now have the following definition:

\begin{definition}[Successor]
	Let $s$ and $s'$ be two periodic orbits.  Let 
	$x\in \F^n_2$ be a state in $s$, and $x' \in \I(x)$.  If the trajectory of the dynamics, with $x'$ the initial condition, enters into $s'$ (in finite time steps), then we say that $s'$ is a successor of $s$. 
\end{definition}

It then naturally leads to the following definition:

\begin{definition}[Stability structure]\label{defH}
	The {\bf stability structure} of a conjunctive Boolean network is  a digraph $H=(S,A)$, with the vertex set being the set of periodic orbits. The edge set of $H$ is defined as follows: Let $s_i$ and $s_j$ be in $S$. Then, $s_is_j$ is an edge of $H$ if  $s_j$ is a successor of $s_i$. Furthermore, an edge $s_is_j$ of $H$ is a {\bf down-edge} (resp. an {\bf up-edge}) if $\sigma(s_i) > \sigma(s_j)$ (resp, $\sigma(s_i) < \sigma(s_j)$).
\end{definition}

Our goal here is to determine the edge set $A$ of $H$. To proceed, we first introduce a {\em partial order} on the set of binary necklaces of length $p^*$: Let $s = y_0\ldots y_{p^*-1} $ and $s' = y'_0\ldots y'_{p^*-1}$ be two binary necklaces. We say that $s$ is {\bf greater than} $s'$, or simply write  $s \succ s'$, if we can obtain $s$ by replacing at least one ``0'' in $s'$ with ``1''. For example, if $s=11100$ and $s'=11000$, then we can obtain $s$ by replacing the third bit ``0'' in $s'$ with ``1'', and thus $s\succ s'$. If, instead, $s=11010$, then there is no way to obtain $s$ by replacing some ``0'' in $s'$ with ``1'', and thus $s$ and $s'$ are not comparable.

With the  definitions and notation above, we state the main result of {\color{black}this} section as follows:

\begin{thm}\label{thm1}
	Let $D$ be the dependency graph associated with a conjunctive Boolean network, and $H=(S,A)$ be the stability structure. Let $s_i$ and $s_j$ be two vertices of $H$. Then,  there is an edge from $s_i$ to $s_j$ if and only if one of the following three conditions is satisfied:  
	\begin{enumerate}
		\item Down-edges: $s_i \succ  s_j$ and $\sigma(s_i)-\sigma(s_j)=1$.  
		\item Up-edges: $s_i \prec s_j$, $\sigma(s_j)-\sigma(s_i)=1$, and $D$ has to be a rose.
		\item Self-loops: $s_i = s_j$, $s_i \neq 1\ldots 1$, and $D$ is not a cycle digraph. 
	\end{enumerate}\,
\end{thm}


We state here a fact as a corollary to Theorem~\ref{thm1}:

\begin{corollary}\label{corthm1}
	Let $D_1$ and $D_2$ be two dependency graphs associated with two conjunctive Boolean networks, having the same loop number $p^*$. Let $H_1$ and  $H_2$ be the corresponding stability structures. Then, $H_1 = H_2$ if one of the following three conditions hold: 
	\begin{enumerate}
		\item Neither $D_1$ nor $D_2$ is a rose.
		\item Both $D_1$ and $D_2$ are roses, but not cycle digraphs. 
		\item Both $D_1$ and $D_2$ are cycle digraphs (and hence $D_1 = D_2$).  
	\end{enumerate}\,
\end{corollary}

\begin{figure}[h]
	\centering
	\includegraphics[height=60mm]{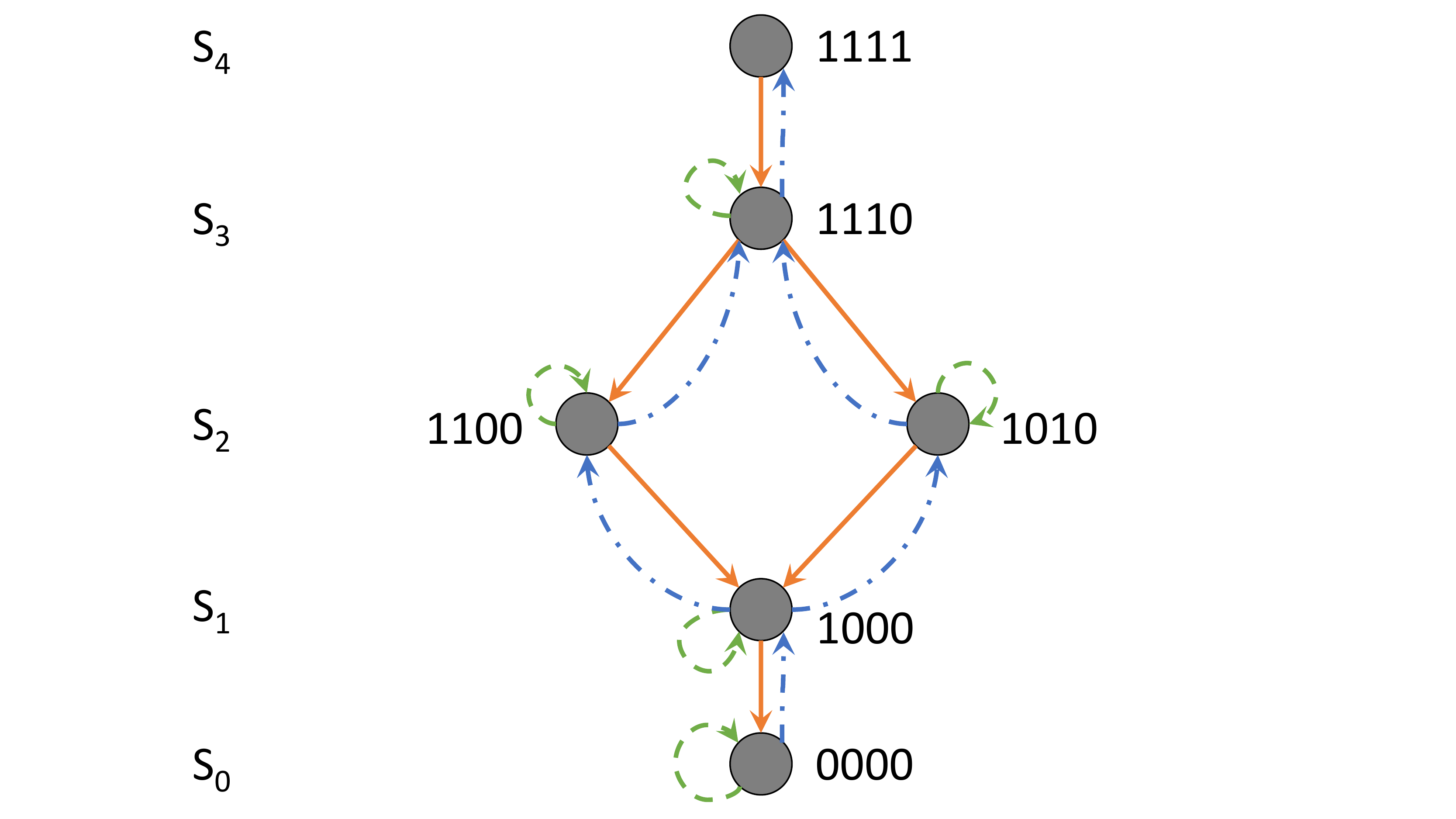}
	\caption{The stability structure $H$ for a dependency graph $D$ with the loop number $p^*=4$. {\color{black}Each edge in $H$ represents a possible transition of the periodic orbit after a single perturbation.} The up-edges exist only for the cases when $D$ is a rose, and the self-loops exist only for the cases when $D$ is not a cycle digraph.}
	\label{H}
\end{figure}

We omit the proof of the corollary as it directly follows from Theorem~\ref{thm1}. We provide an example in Fig. \ref{H} for the case when the loop number $p^* = 4$.

{\color{black}
\begin{Remark}\label{complexity}
For a strongly connected graph which is neither a rose nor a cycle digraph, the loop number uniquely determines the periodic orbits and the stability structure. 
Thus, computing its stability structure is reduced to computing the loop number, which is the greatest common divisor of lengths of all cycles. We note here a few research works on finding the cycles of an arbitrary digraph~\cite{johnson1975finding,weinblatt1972new,mateti1976algorithms,bax1994algorithms}. For example,~\cite{johnson1975finding} proposed an algorithm which finds all cycles of a digraph in time bounded by $\mathcal{O}((n+e)(N+1))$, where $n$ is the number of vertices, $e$ is the number of edges, and $N$ is the number of cycles in the digraph. There are also algorithms for finding the cycles of specific lengths in a digraph~\cite{bax1996finite,alon1997finding,czumaj2014finding,giscard2016general}. Connecting these computational complexity results  to the structural results in this paper would be a fruitful direction of future research, as also discussed in the Conclusions section. 
\end{Remark}
}

The remainder of this subsection is devoted to the proof of Theorem~\ref{thm1}. 

Let $s = y_0\ldots y_{p^*-1}$ be a periodic orbit, and $x$ be a state in $s$. Let $x' \in \I(x)$, with $x'_1 \neq x_1$. Let $G_k = (U_k, F_k)$, for $k = 0,\ldots, p^*-1$, be irreducible components of $D$.  Then, from Corollary~\ref{same}, we can assume without loss of generality that 
$$
x_{U_k} = y_k {\bf 1}, \hspace{10pt} \forall k = 0,\ldots, p^*-1,
$$
where ${\bf 1}$ is a vector of all ones with an appropriate dimension. We may further assume that $x_1$ is an entry of $x_{U_0}$. So, $x_1 = y_0$ and $x'_1 = \neg y_0$ (negating the value of $y_0$). 
With these preliminaries, we establish the following result:

\begin{proposition}\label{pro:forthm1}
	Let $s$, $x$ and $x'$ be defined as above. Suppose that the trajectory, with $x'$ the initial condition, {\color{black}enters the periodic orbit $s'$}. Then, there are two cases:  
	\begin{enumerate}
		\item If $|U_0| = 1$, then $s' = (\neg y_0)y_1\ldots y_{p^*-1}$.
		\item If $|U_0| > 1$, then $s' = 0 y_1\ldots y_{p^*-1}$.
	\end{enumerate}\,
\end{proposition}

\noindent
{\bf Proof.}
	For the case $|U_0| = 1$, we note from Corollary~\ref{same} that the state $x'$ is already in the periodic orbit $s'$. We now prove for the case $|U_0| > 1$. 
	
	First, note that the vector $x'_{U_0}$, obtained by restricting $x'$ to $U_0$, must contain an entry of value $0$. This holds because if $y_0 = 0$, then from Corollary~\ref{same}, all the entries  of $x_{U_0}$ hold value $0$. Since $x'_{U_0}$ is derived by negating the value of $x_1$,  there are $(|U_0| - 1)$ zeros in $x'_{U_0}$. If $y_0 = 1$, then by construction, $x'_1 = 0$, which is contained in $x'_{U_0}$.   
	
	Next, consider the induced dynamics on $G_0$: First, from the value update rule and the first item of Theorem~\ref{dynamic},  if $x'_{U_0}(0)$ contains an entry of value $0$, then so does $x'_{U_0}(tp^*)$ for all $t \ge 0$. Second, since $G_0$ is irreducible, a periodic orbit of the induced dynamics has to be a fixed point. Combining these two facts, we know that there is a time $t_0\ge 0$ such that $x'_{U_0}(tp^*) = {\bf 0}$ for all $t \ge t_0$.
	
	Now, for each $k = 1,\ldots, p^*-1$, we appeal again to the first item of Theorem~\ref{dynamic} and obtain
	$$
	x'_{U_k}(tp^*) = f^{tp^*}_{U_k}(x'_{U_k}(0)) = g^t_k(x'_{U_k}(0)) = x'_{U_k}(0).
	$$
	The last equality holds because by construction of $x'$, we have that $x'_{U_k}(0) = y_k {\bf 1}$ which is a fixed point of the induced dynamics on $G_k$. The relation above holds for all $t \ge 0$.
	
	Combining the arguments above, we conclude that for any $t \ge t_0$, we have
	$$
	x'_{U_k}(tp^*) =
	\left\{
	\begin{array}{ll}
	0 & \mbox{if } k = 0 \\
	y_k {\bf 1} & \mbox{otherwise}.
	\end{array}
	\right. 
	$$
	Thus, $s' = 0y_1\ldots y_{p^*-1}$.
	\hfill$\qed$

With Proposition~\ref{pro:forthm1}, we prove Theorem~\ref{thm1}:

\noindent
{\bf Proof of Theorem~\ref{thm1}.}
	Let $s_i = y_0\ldots y_{p^*-1}$ and $s_j = y'_0\ldots y'_{p^*-1}$. There are three cases to consider: 
	
	{\em Case 1: $\sigma(s_i) > \sigma(s_j)$.}  If $s_is_j$ is an edge of $H$, then from the proof of Proposition~\ref{pro:forthm1},  we must have $y_0 = 1$, $y'_0 = 0$, and $y_i = y'_i$ for all $i =1,\ldots, p^*-1$ (after appropriate rotations of the strings). In other words, we have 
	$s_i \succ s_j$ and  $\sigma(s_i) - \sigma(s_j) = 1$.  Conversely, if the condition in the first item of the theorem is satisfied, then we can always write $s_i = 1y_1\ldots y_{p^*-1}$ and $s_j = 0y_1\ldots y_{p^*-1}$.  Let $x$ be a state in $s_i$, with $x_{U_0} = {\bf 1}$ and $x_{U_k} = y_k {\bf 1}$ for all $k = 1,\ldots, p^*-1$. Then, by negating the value of an entry of $x_{U_1}$, we obtain a state $x'$. Moreover, from Proposition~\ref{pro:forthm1}, the trajectory, with $x'$ the initial condition, will enter $s_j$ in finite time steps, and hence $s_j$ is a successor of $s_i$.

	{\em Case 2: $\sigma(s_i) < \sigma(s_j)$.} If $s_is_j$ is an edge of $H$, then from Proposition~\ref{pro:forthm1}, we must have $s_i \prec s_j$ and $\sigma(s_j) - \sigma(s_i) = 1$, and moreover there exists at least one $k$ such that $|U_k| =1$. Then, from the first item of Proposition~\ref{special}, $D$ has to be a rose. Conversely, if the condition in the second item of the theorem is satisfied, then again, by the first item of Proposition~\ref{special}, there is a $k$ such that $|U_k| = 1$. Without loss of generality, we assume that $k = 0$, and write $s_i = 0y_1\ldots y_{p^*-1}$ and $s_j = 1y_1\ldots y_{p^*-1}$. 
	Let $x$ be a state in $s_i$ with $x_{U_0} = 0$ (note that $x_{U_0}$ is a scalar in this case). By negating the value of $x_{U_0}$, we obtain a new state $x'$. Then, from Proposition~\ref{pro:forthm1}, the trajectory, with $x'$ the initial condition, will enter the periodic orbit $s_j$. Thus, $s_j$ is a successor of $s_i$.

	{\em Case 3: $\sigma(s_i) = \sigma(s_j)$.} If $s_is_j$ is an edge of $H$, then from Proposition~\ref{pro:forthm1}, we must have that (i) $s_i = s_j$; (ii) there exists at least one $k = 0,\ldots, p^*-1$ such that $y_k = 0$ and $|U_k| > 1$. Combining the condition (ii) and the second item of Proposition~\ref{special}, we know that $D$ can not be a cycle digraph and $s_i \neq 1\ldots 1$. Conversely, if the condition in the third item of the theorem is satisfied, then there is a $k$ such that $|U_k| > 1$. Without loss of generality, we assume that $k = 0$, and hence $s_i $ can be written as $0 y_1\ldots y_{p^*-1}$.  Since $s_i \neq 1\ldots 1$, we can find a state $x$ in the periodic orbit $s_i$ such that $x_{U_0} = {\bf 0}$. Now, by negating the value of an entry of $x_{U_0}$, we obtain a new state $x'$. Then, from Proposition~\ref{pro:forthm1}, the trajectory, with $x'$ the initial condition, will enter the periodic orbit $0y_1\ldots y_{p^*-1}$, which is $s_i$ itself.
	\hfill$\qed$   


\subsection{Transition weights}\label{transitionweights}
In this subsection, we introduce and compute the transition weight for each edge of the stability structure $H$. First, recall that the set $\I(x)$ is comprised of the states that differ from $x$ by only one entry. It should be clear that $|\I(x)| = n$ for all $x\in \F^n_2$. Now, let $s_i = \{x(t_0), \ldots, x(t_0 + p-1)\}$ be a periodic orbit, and $s_j$ be a  successor of $s_i$. We define $\mu(s_i,s_j)$ to be the total number of pairs $(x, x')$, for $x\in s_i$ and $x' \in \I(x)$, such that the trajectory of the conjunctive Boolean network,  with $x'$ the initial condition, enters into $s_j$. We then have the following definition:

\begin{definition}[Transition weight]
	Let $s_i$ be a periodic orbit of period $p$, and $s_j$ be its successor. Then, the transition weight $P(s_i,s_j)$ on the edge $s_is_j$ of the stability structure $H$ is $$P(s_i,s_j):= \frac{\mu(s_i,s_j)}{np}.$$\,
\end{definition}

We note here that by the definition, $\sum_{s_j} P(s_i,s_j) = 1$, where the summation is over the successors of $s_i$. Thus, each $P(s_i,s_j)$ can be understood as the probability of the transition from $s_i$ to $s_j$ upon the condition that the pair $(x,x')$ is uniformly chosen from  the set $\{(x,x') \mid x\in s_i, x'\in \I(x)\}$.

\begin{figure}[h]
	\centering
	\includegraphics[height=60mm]{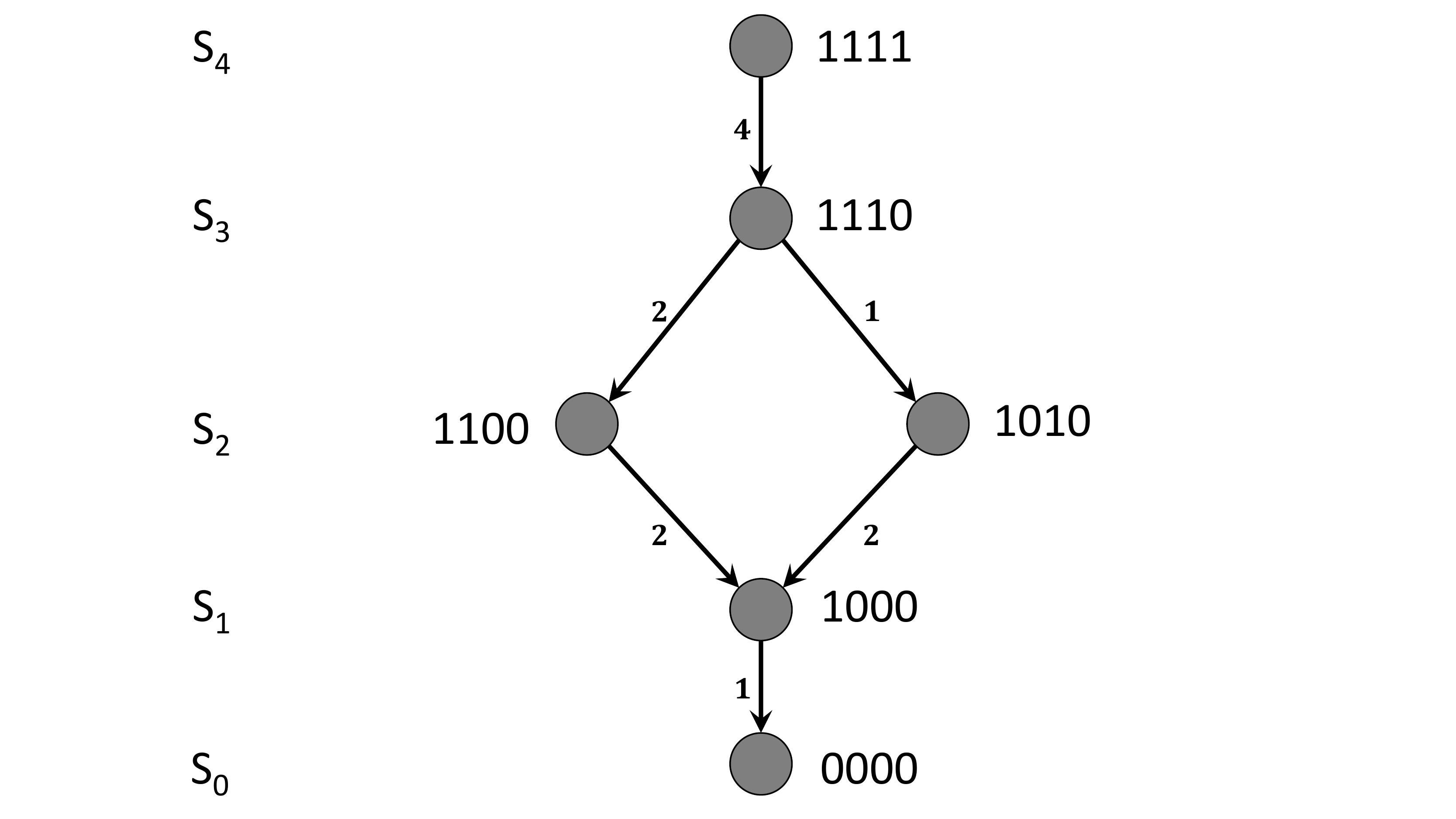}
	\caption{The values of $\gamma$'s for a dependency graph $D$ with the loop number $p^*=4$. The number labeled on an $ss'$ edge is the value of $\gamma(s,s')$.}
	\label{gamma}
\end{figure}

For the remainder of this subsection, we evaluate the transition weight $P(s_i,s_j)$. 
To proceed, first note that by the arguments in the beginning of Subsection~\ref{Identify}, we can identify the two periodic orbits $s_i$ and $s_j$ with two binary necklaces: $s_i = y_0\ldots y_{p^*-1}$ and $s_j = y'_0\ldots y'_{p^*-1}$. From Theorem~\ref{thm1}, we know that one of the following three conditions holds:
\begin{enumerate}
	\item $s_i \succ s_j$ and $\sigma(s_i) - \sigma(s_j) = 1$. 
	\item $s_i = s_j$, $\sigma(s_i)\neq p^*$ and $D$ is not a cycle digraph.
	\item $s_i \prec s_j$, $\sigma(s_j) - \sigma(s_i) = 1$ and $D$ is a rose.
\end{enumerate}
We thus introduce the following number for a pair of necklaces: 
Let $s,s'$ be two necklaces of equal length $p^*$, with $s \succ s'$ and $\sigma(s)-\sigma(s')=1$. We define $\gamma(s,s')$ to be the number of ways to obtain $s'$ from $s$ by replacing a ``1'' in $s$ with a ``$0$''. We note here that from its definition, $\gamma(s,s')$ can also be viewed as the number of ways to obtain $s$ from $s'$ by replacing a ``0'' in $s'$ with a ``1''. For example, consider the case where $p^* = 4$, and $s=1110$, $s'=1100$. Then, there are two ways to obtain $s'$ from $s$: One way is to replace the first ``$1$'' in $s$ with ``$0$''. The other way is to  replace the third ``$1$'' with ``$0$''. So, in this case, $\gamma(s,s')  = 2$. 
We also refer to Fig.~\ref{gamma} for other values of $\gamma(s,s')$ under the case $p^* = 4$. We further note that $\sum_{s'} \gamma(s,s') = \sigma(s)$, where the summation is over the successors of $s$ other than itself. 


We further need the following definition: Given a digraph $D$. We define an integer number $\alpha$ as follows: We let $\alpha= 0$ if $D$ is not a rose. Otherwise, we let $\alpha$ be the number of common vertices of the cycles of the rose $D$. From the proof of Proposition~\ref{special}, we know that $\alpha$ is also the number of irreducible components comprised of a single vertex.

With the definitions and notation above, we establish the following result.

\begin{proposition}\label{thm2}
	Let $D$ be a dependency graph associated with a conjunctive Boolean network. 
	Let $s_i$ and $s_j$ be two periodic orbits (which can be identified as two binary necklaces). Then, the following holds: 
	\[
	P(s_i,s_j) =
	\begin{cases}{}
	
	\vspace{0.1cm}
	
	\frac{\gamma(s_i,s_j)}{p^*} &
	\begin{array}{lr}
	\text{if } s_i\succ s_j\\ 
	\text{and } \sigma(s_i)-\sigma(s_j)=1,
	\end{array}\\

	\vspace{0.1cm}
	
	\frac{\alpha\gamma(s_j,s_i)}{np^*} & 
	\begin{array}{lr}
	\text{if } s_i\prec s_j\\ 
	\text{and } \sigma(s_j)-\sigma(s_i)=1,
	\end{array}\\
	
	\vspace{0.1cm}
	
	\frac{(p^*-\sigma(s_i))(n-\alpha)}{np^*} & 
	\begin{array}{lr}
	\text{if } s_i= s_j\\ 
	\text{and } \sigma(s_i)\neq p^*,
	\end{array}\\

	\vspace{0.1cm}

	0 & 
	\begin{array}{lr}
	\text{otherwise}.
	\end{array}
	\end{cases}
	\]\,
\end{proposition}

\begin{figure}[h]
	\centering
	\includegraphics[height=60mm]{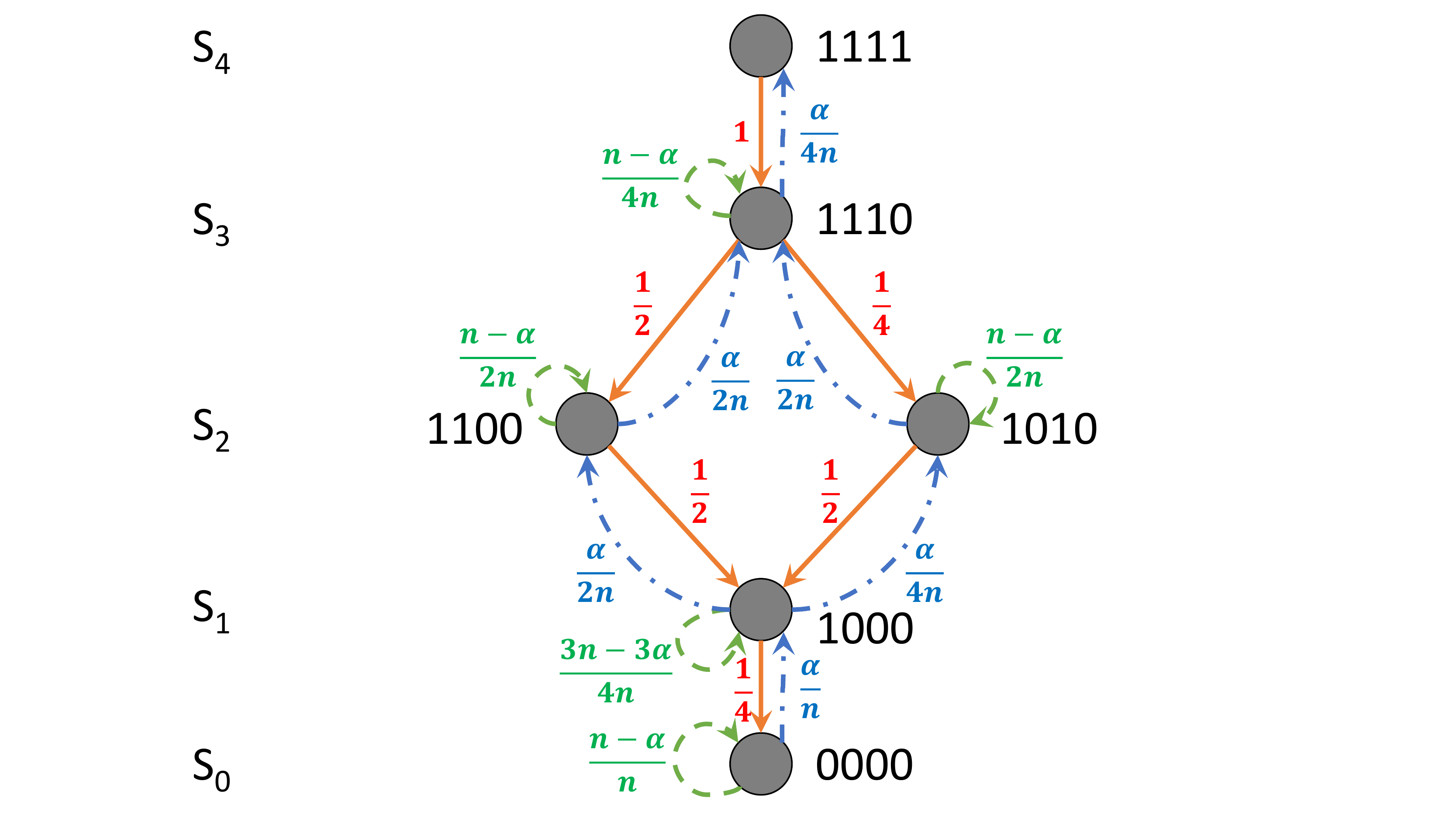}
	\caption{The values of the transition weights $P(s_i,s_j)$ are labeled on the edge of the stability structure $H$ for the case where $D$ has loop number $p^*=4$.}
	\label{H_prob}
\end{figure}

\begin{Remark}
	We note here that if the graph $D$ is not a rose, then the transition weights depend only on the loop number of $D$. We provide an example in Fig.~\ref{H_prob} for the case  $p^*=4$.
\end{Remark}

\noindent
{\bf Proof of Proposition~\ref{thm2}.}
	Let $s_i = \{x(t_0),\ldots,x(t_0 + p-1)\}$ be a periodic orbit, and $s_j$ be a  successor of $s_i$.
	We first introduce some notations. For a time step $t$, for $t_0 \le t \le t_0 + p^*-1$, we let $\mu^*_t(s_i,s_j)$ be the number of pairs $(x(t), x')$,  for $x'\in \I(x(t))$, such that the trajectory of the conjunctive Boolean network,  with $x'$ the initial condition, enters into $s_j$. We further let $\mu^*(s_i,s_j) = \sum^{t_0 + p^*-1}_{t = t_0}\mu^*_t(s_i,s_j)$. Then, from its definition, we have the following relation: 
	$$
	P(s_i,s_j)=\frac{\mu(s_i,s_j)}{np}=\frac{\mu^*(s_i,s_j)}{np^*}.
	$$
	We now evaluate $\mu_t^*(s_i,s_j)$ for each $t = t_0,\ldots, t_0 + p^*-1$.

	Let $G_k = (U_k, F_k)$, for $k = 0,\ldots, p^*-1$, be irreducible components of $D$. From Corollary~\ref{same}, we can write $x_{U_k}(t) = y_k(t) {\bf 1}$ for all $k = 0,\ldots, p^*-1$. We identify $s_i$ with $y_0(t)\ldots y_{p^*-1}(t)$, and $s_j$ with $y'_0\ldots y'_{p^*-1}$. 
	We further let $\Gamma(t)$ be a subset of $\{0,\ldots, p^*-1\}$ defined as follows: an index $k$ is in $\Gamma(t)$ if the binary necklace $y'_0\ldots y'_{p^*-1}$ can be obtained from $y_0(t)\ldots y_{p^*-1}(t)$ by negating the value of $y_k(t)$. 
	
	We now relate $\Gamma(t)$ and $\Gamma(t')$ for two different time steps $t$ and $t'$. In particular, we show that if $k\in \Gamma(t)$, then 
	\begin{equation}\label{relation99}
	\left ((k + t' - t) \bmod p^*\right ) \in \Gamma(t').
	\end{equation}
	To see this, note that from the second item of Theorem~\ref{dynamic}, we have that for all $k = 0,\ldots, p^*-1$, 
	$x_{U_{(k+t' -t) \bmod p^*}}(t') = x_{U_k}(t),$
	and hence 
	$$y_{(k+t' -t) \bmod p^*}(t') = y_k(t),$$
	which implies~\eqref{relation99}.
	
	We first evaluate the transition weights for down-edges. Note that since $s_i \succ s_j$, $y_k(t) = 1$ for all $k\in \Gamma(t)$, and moreover, $|\Gamma(t)| = \gamma(s_i,s_j)$.             
	To proceed, we fix an index $k \in \{0,\ldots, p^*-1\}$, and assume that $x_{U_k}(t) = {\bf 1}$, and let $x' \in \I(x(t))$ be derived by negating an entry of $x_{U_k}$. Then, from Proposition~\ref{pro:forthm1}, the trajectory, with $x'$ the initial condition, will enter a periodic orbit, which can be identified as the following binary necklace: $y_0(t) \ldots y_{k-1}(t)\,0 \,y_{k+1}(t) \ldots y_{p^*-1}(t).$ So, $y'_0\ldots y'_{p^*-1}$ coincides with the binary necklace above if and only if $k \in \Gamma(t)$. It then follows that
	\begin{equation}\label{relation88}
	\mu^*_t(s_i,s_j) = \sum_{k\in \Gamma(t)} |U_k|.
	\end{equation}     
	
	Then, by combining~\eqref{relation99} and~\eqref{relation88} with the fact that $|\Gamma(t)| = \gamma(s_i,s_j)$ for all~$t$, we obtain
	$$
	\begin{array}{lll}
	\mu^*(s_i,s_j) & = & \sum^{p^*-1}_{t=0} \mu_t^*(s_i,s_j)   \vspace{3pt} \\
	& = & \sum^{p^*-1}_{t = 0}\sum_{k\in \Gamma(t)} |U_k|  \vspace{3pt} \\
	& = & \gamma(s_i,s_j) \sum^{p^*-1}_{k = 0} |U_k| = \gamma(s_i,s_j) n.
	\end{array}
	$$
	Thus, for the case $s_i\succ s_j$ and $\sigma(s_i) - \sigma(s_j) = 1$, we obtain $P(s_i,s_j) = \gamma(s_i,s_j)/p^*$.

	We next evaluate the transition weights for up-edges. Note that since $s_i\prec s_j$, $y_k(t)=0$ for all $k\in\Gamma(t)$, and moreover, $|\Gamma(t)|=\gamma(s_j,s_i)$. To proceed, we fix an index $k\in \{0,\ldots,p^*-1\}$, and assume that $x_{U_k}(t)={\bf 0}$, and let $x' \in \I(x(t))$ be derived by negating an entry of $x_{U_k}$. We know from Proposition~\ref{pro:forthm1}  that if $|U_k| = 1$,  then the trajectory, with $x'$ the initial condition, will enter a periodic orbit identified as: $y_0(t) \ldots y_{k-1}(t)\,1 \,y_{k+1}(t) \ldots y_{p^*-1}(t),$ and vice versa. So, $y'_0\ldots y'_{p^*-1}$ coincides with the binary necklace above if and only if $k \in \Gamma(t)$ and $|U_k|=1$. It then follows that
	\begin{equation}\label{relation77}
	\mu^*_t(s_i,s_j) = \sum_{k:\left\{\substack{k\in \Gamma(t)\\|U_k|=1}\right\}} |U_k|.
	\end{equation}     
	
	Then, by combining~\eqref{relation99} and~\eqref{relation77} with the fact that $|\Gamma(t)| = \gamma(s_j,s_i)$ for all~$t$, we obtain
	$$
	\begin{array}{lll}
	\mu^*(s_i,s_j) & = & \sum^{p^*-1}_{t=0} \mu_t^*(s_i,s_j)   \vspace{3pt} \\
	& = & \sum^{p^*-1}_{t = 0}\sum_{k:\left\{\substack{k\in \Gamma(t)\\|U_k|=1}\right\}} |U_k|  \vspace{3pt} \\
	& = & \gamma(s_j,s_i) \sum_{\substack{k=0,\ldots, p^*-1\\|U_k|=1}} |U_k| = \gamma(s_j,s_i) \alpha.
	\end{array}
	$$
	Thus, for the case $s_i\succ s_j$ and $\sigma(s_i) - \sigma(s_j) = 1$, we obtain $P(s_i,s_j) = \alpha\gamma(s_j,s_i)/(np^*)$. 
	
	As a final step, we evaluate the transition weights for self-loops. 
	To proceed, we fix an index $k\in \{0,\ldots,p^*-1\}$, and assume that $x_{U_k}(t)={\bf 0}$, and let $x' \in \I(x(t))$ be derived by negating an entry of $x_{U_k}$. We know from Proposition~\ref{pro:forthm1} that if $|U_k| > 1$, then the trajectory, with $x'$ the initial condition, will enter a periodic orbit identified as: $y_0(t) \ldots y_{k-1}(t)\,0 \,y_{k+1}(t) \ldots y_{p^*-1}(t),$ and vice versa. So, $y'_0\ldots y'_{p^*-1}$ coincides with the binary necklace above if and only if $x_{U_k}(t)={\bf 0}$ and $|U_k|>1$. It then follows that
	\begin{equation}\label{relation66}
	\mu^*_t(s_i,s_j) = \sum_{k:\left\{\substack{x_{U_k}(t)={\bf 0}\\|U_k|>1}\right\}} |U_k|.
	\end{equation}
	Following  the second item of Theorem~\ref{dynamic},  we have that if $k\in \{k:x_{U_k}(t)={\bf 0}\}$, then
	\begin{equation}\label{relation55}
	\left ((k + t' - t) \bmod p^*\right ) \in \{k:x_{U_k}(t')={\bf 0}\}.
	\end{equation}
	Now, by combining~\eqref{relation66} and~\eqref{relation55} with the fact that $|\{k:x_{U_k}(t)={\bf 0}\}| = (p^*-\sigma(s_i))$ for all~$t$, we obtain
	$$
	\begin{array}{lll}
	\mu^*(s_i,s_j) & = & \sum^{p^*-1}_{t=0} \mu_t^*(s_i,s_j)   \vspace{3pt} \\
	& = & \sum^{p^*-1}_{t = 0}\sum_{k:\left\{\substack{x_{U_k}(t)={\bf 0}\\|U_k|>1}\right\}} |U_k|  \vspace{3pt} \\
	& = & (p^*-\sigma(s_i)) \sum_{\substack{k=0,\ldots, p^*-1\\|U_k|> 1}} |U_k| \vspace{3pt}\\
	& = & (p^*-\sigma(s_i))(n-\alpha).
	\end{array}
	$$
	Thus, for the case $s_i=s_j$ and $\sigma(s_i)\neq p^*$, we obtain $P(s_i,s_j) = (p^*-\sigma(s_i))(n-\alpha)/(np^*)$. 
	\hfill$\qed$


\section{Conclusions and Outlooks} \label{end}

In this paper, we have introduced a new approach to study the dynamics of a conjunctive Boolean network, and investigated the stability structure of the periodic orbits. Specifically, we have proposed a vertex set partition of the dependency graph associated with the conjunctive Boolean network, and decomposed the digraph into multiple irreducible components. We have then introduced the induced dynamics on each of the irreducible components, and established in Theorem~\ref{dynamic} the relationship  between the original conjunctive Boolean network and the induced dynamics on the irreducible components. Following this relationship, we have further identified the periodic orbits of the conjunctive Boolean network with binary necklaces. By introducing a partial ordering on the set of binary necklaces, we have established in Theorem~\ref{thm1} the stability structure of the periodic orbits. In particular, we have provided in the theorem a  necessary and sufficient condition for the existence of a transition from one periodic orbit to another under the condition that a single perturbation occurs to a state of the periodic orbit. The transition weights are also evaluated in Subsection~\ref{transitionweights}. 

Although systems biology serves as the main motivation for our research, applications of this work are by far not limited to {\color{black}gene regulatory networks}. Conjunctive Boolean networks are also suitable to model, for example, water quality networks. In such networks, each Boolean variable can be viewed as the water quality within a pipe. The Boolean variable takes ``1'' if the water is not polluted, and ``0'' if the water is polluted. The water in each pipe comes from some other pipes, and is polluted if the water in any one of those other pipes was polluted. Other examples which can be modeled by conjunctive Boolean networks include social networks (information flow on Twitter or Facebook) and supply chain networks (movement of materials),  and the results of this paper can be applied to all these networks as well.


{\color{black} There are a number of research directions in this area that could be pursued in the future, and we mention here a few of them. First, we recall that the stability structure is constructed by assuming that only a single entry of a state in a periodic orbit is perturbed. A natural question is then to ask what the stability structure would be like if more than one entry of the state are perturbed. We will investigate this question in the future. {\color{black}We also recall from Remark~\ref{complexity} that the loop number determines the set of periodic orbits, and also the stability structure (provided that the dependency graph is not a rose). We also provided a few references on algorithms of finding all cycles of a general digraph. A future direction following this would be to develop/improve algorithms for computing the loop number (so as to construct the stability structure) of strongly connected graphs, and to evaluate the computational complexity of these algorithms. }
	
	We also note here that conjunctive Boolean networks have been studied mostly  over strongly connected digraphs. It is still not clear how to identify the periodic orbits of a weakly connected conjunctive Boolean network. One promising approach is to apply the strong component decomposition to the weakly connected digraph {\color{black}(see, for example,~\cite{chen2015controllability,jarrah2010dynamics})} which partitions the digraph into strongly connected subgraphs. A few results obtained in this paper  can be used to establish certain properties of a periodic orbit when restricted to each connected component. Yet, a complete understanding of the asymptotic behavior is  still lacking, let alone characterizing the stability structure. 
}


Amongst other future directions, we mention the research on the dynamics of non-conjunctive Boolean networks. For example, we have been working on Boolean networks whose value update rules are given by XOR (XNOR) and NAND (NOR) operations. Note, in particular, that  the approach proposed in this paper---decomposing the dependency graph into irreducible components and relating the original dynamic to induced dynamics---can still be applied to any of these Boolean networks. 

Last but not least, we are also interested in the control of a Boolean network. Consider, for example, that there is a subset of vertices whose values can be manipulated at any time step while the other vertices follow a certain value update rule. A question of particular interest is whether the system is controllable or not, i.e., for any initial condition $x_0$ and any final condition $x_1$, whether there is a time duration $T$, and control sequences for the manipulable vertices such that the Boolean network can be steered from $x_0$ to $x_1$ after $T$ time steps. {\color{black} This question has been addressed in~\cite{orbitcontrol,statecontrol,gaocontrollability}.} Other questions, such as minimal controllability, i.e., finding a subset of vertices with least cardinality such that the system is controllable by manipulating the values of these vertices, are also in the scope of our future work.






\bibliographystyle{plain}        
\bibliography{autosam}           

\begin{thebibliography}{10}

\bibitem{alcolei2015flora}
A.~Alcolei, K.~Perrot, and S.~Sen{\'e}.
\newblock On the flora of asynchronous locally non-monotonic {B}oolean automata
  networks.
\newblock {\em arXiv preprint arXiv:1510.05452}, 2015.

\bibitem{alon1997finding}
N.~Alon, R.~Yuster, and U.~Zwick.
\newblock Finding and counting given length cycles.
\newblock {\em Algorithmica}, 17(3):209--223, 1997.

\bibitem{bax1996finite}
E.~Bax and J.~Franklin.
\newblock A finite-difference sieve to count paths and cycles by length.
\newblock {\em Information Processing Letters}, 60(4):171--176, 1996.

\bibitem{bax1994algorithms}
E.~T. Bax.
\newblock Algorithms to count paths and cycles.
\newblock {\em Information Processing Letters}, 52(5):249--252, 1994.

\bibitem{chen2015consensus}
X.~Chen, M-A. Belabbas, and T.~Ba\c{s}ar.
\newblock Consensus with linear objective maps.
\newblock In {\em Proc. 54th Conference on Decision and Control (CDC)}, pages
  2847--2852. IEEE, 2015.

\bibitem{chen2015controllability}
X.~Chen, M.-A. Belabbas, and T.~Ba{\c s}ar.
\newblock Controllability of formations over directed time-varying graphs.
\newblock {\em IEEE Transactions on Control of Network Systems}, 2015.

\bibitem{colon2006monomial}
O.~Col{\'o}n-Reyes, A.~Jarrah, R.~Laubenbacher, and B.~Sturmfels.
\newblock Monomial dynamical systems over finite fields.
\newblock {\em arXiv preprint math/0605439}, 2006.

\bibitem{colon2005boolean}
O.~Col{\'o}n-Reyes, R.~Laubenbacher, and B.~Pareigis.
\newblock Boolean monomial dynamical systems.
\newblock {\em Annals of Combinatorics}, 8(4):425--439, 2005.

\bibitem{czumaj2014finding}
A.~Czumaj, O.~Goldreich, D.~Ron, C.~Seshadhri, A.~Shapira, and C.~Sohler.
\newblock Finding cycles and trees in sublinear time.
\newblock {\em Random Structures \& Algorithms}, 45(2):139--184, 2014.

\bibitem{etesami2016complexity}
S.~R. Etesami and T.~Ba{\c{s}}ar.
\newblock Complexity of equilibrium in competitive diffusion games on social
  networks.
\newblock {\em Automatica}, 68:100--110, 2016.

\bibitem{funahashi1993approximation}
K.~Funahashi and Y.~Nakamura.
\newblock Approximation of dynamical systems by continuous time recurrent
  neural networks.
\newblock {\em Neural networks}, 6(6):801--806, 1993.

\bibitem{gaocontrollability}
Z.~Gao, X.~Chen, and T.~Ba{\c{s}}ar.
\newblock Controllability of conjunctive {B}oolean networks with application to
  gene regulation.
\newblock 2017.
\newblock submitted to \emph{IEEE Transactions on Control of Network Systems}.

\bibitem{orbitcontrol}
Z.~Gao, X.~Chen, and T.~Ba{\c{s}}ar.
\newblock Orbit-controlling sets for conjunctive {B}oolean networks.
\newblock In {\em Proc. 2017 American Control Conference (ACC)}, pages
  4989--4994, 2017.

\bibitem{statecontrol}
Z.~Gao, X.~Chen, and T.~Ba{\c{s}}ar.
\newblock State-controlling sets for conjunctive {B}oolean networks.
\newblock In {\em Proc. 20th IFAC World Congress}, pages 14855--14860, 2017.

\bibitem{gao2016periodic}
Z.~Gao, X.~Chen, J.~Liu, and T.~Ba{\c{s}}ar.
\newblock Periodic behavior of a diffusion model over directed graphs.
\newblock In {\em Proc. 55th Conference on Decision and Control (CDC)}, pages
  37--42. IEEE, 2016.

\bibitem{georgescu2008gene}
C.~Georgescu, W.~Longabaugh, D.~D. Scripture-Adams, E.~David-Fung, M.~A. Yui,
  M.~A. Zarnegar, H.~Bolouri, and E.~V. Rothenberg.
\newblock A gene regulatory network armature for {T} lymphocyte specification.
\newblock {\em Proceedings of the National Academy of Sciences},
  105(51):20100--20105, 2008.

\bibitem{gilbert1961symmetry}
E.~N. Gilbert and J.~Riordan.
\newblock Symmetry types of periodic sequences.
\newblock {\em Illinois Journal of Mathematics}, 5(4):657--665, 1961.

\bibitem{giscard2016general}
P.-L. Giscard, N.~Kriege, and R.~C. Wilson.
\newblock A general purpose algorithm for counting simple cycles and simple
  paths of any length.
\newblock {\em arXiv preprint arXiv:1612.05531}, 2016.

\bibitem{goles2012disjunctive}
E.~Goles and M.~Noual.
\newblock Disjunctive networks and update schedules.
\newblock {\em Advances in Applied Mathematics}, 48(5):646--662, 2012.

\bibitem{harris2002model}
S.~E. Harris, B.~K. Sawhill, A.~Wuensche, and S.~Kauffman.
\newblock A model of transcriptional regulatory networks based on biases in the
  observed regulation rules.
\newblock {\em Complexity}, 7(4):23--40, 2002.

\bibitem{imer2006optimal}
O.~C. Imer, S.~Y{\"u}ksel, and T.~Ba{\c{s}}ar.
\newblock Optimal control of {LTI} systems over unreliable communication links.
\newblock {\em Automatica}, 42(9):1429--1439, 2006.

\bibitem{jarrah2010dynamics}
A.~S. Jarrah, R.~Laubenbacher, and A.~Veliz-Cuba.
\newblock The dynamics of conjunctive and disjunctive {B}oolean network models.
\newblock {\em Bulletin of Mathematical Biology}, 72(6):1425--1447, 2010.

\bibitem{jarrah2007nested}
A.~S. Jarrah, B.~Raposa, and R.~Laubenbacher.
\newblock Nested canalyzing, unate cascade, and polynomial functions.
\newblock {\em Physica D: Nonlinear Phenomena}, 233(2):167--174, 2007.

\bibitem{johnson1975finding}
D.~B. Johnson.
\newblock Finding all the elementary circuits of a directed graph.
\newblock {\em SIAM Journal on Computing}, 4(1):77--84, 1975.

\bibitem{kauffman1969homeostasis}
S.~Kauffman.
\newblock Homeostasis and differentiation in random genetic control networks.
\newblock {\em Nature}, 224:177--178, 1969.

\bibitem{kauffman2003random}
S.~Kauffman, C.~Peterson, B.~Samuelsson, and C.~Troein.
\newblock Random {B}oolean network models and the yeast transcriptional
  network.
\newblock {\em Proceedings of the National Academy of Sciences},
  100(25):14796--14799, 2003.

\bibitem{kauffman2004genetic}
S.~Kauffman, C.~Peterson, B.~Samuelsson, and C.~Troein.
\newblock Genetic networks with canalyzing {B}oolean rules are always stable.
\newblock {\em Proceedings of the National Academy of Sciences},
  101(49):17102--17107, 2004.

\bibitem{kauffman1969metabolic}
S.~A. Kauffman.
\newblock Metabolic stability and epigenesis in randomly constructed genetic
  nets.
\newblock {\em Journal of Theoretical Biology}, 22(3):437--467, 1969.

\bibitem{stuart1993origins}
S.~A. Kauffman.
\newblock {\em The {O}rigins of {O}rder: Self {O}rganization and {S}election in
  {E}volution}.
\newblock Oxford University Press, USA, 1993.

\bibitem{khanafer2014stability}
A.~Khanafer, T.~Ba\c{s}ar, and B.~Gharesifard.
\newblock Stability properties of infected networks with low curing rates.
\newblock In {\em American Control Conference (ACC), 2014}, pages 3579--3584.
  IEEE, 2014.

\bibitem{mateti1976algorithms}
P.~Mateti and N.~Deo.
\newblock On algorithms for enumerating all circuits of a graph.
\newblock {\em SIAM Journal on Computing}, 5(1):90--99, 1976.

\bibitem{melliti2013convergence}
T.~Melliti, D.~Regnault, A.~Richard, and S.~Sen{\'e}.
\newblock On the convergence of {B}oolean automata networks without negative
  cycles.
\newblock In {\em Cellular Automata and Discrete Complex Systems}, pages
  124--138. Springer, 2013.

\bibitem{mendoza1999genetic}
L.~Mendoza, D.~Thieffry, and E.~R. Alvarez-Buylla.
\newblock Genetic control of flower morphogenesis in {A}rabidopsis {T}haliana:
  a logical analysis.
\newblock {\em Bioinformatics}, 15(7):593--606, 1999.

\bibitem{moreau1872permutations}
C~Moreau.
\newblock Sur les permutations circulaires distinctes.
\newblock {\em Nouvelles Annales de Math{\'e}matiques, Journal des Candidats
  aux {\'E}coles Polytechnique et Normale}, 11:309--314, 1872.

\bibitem{noual2012updating}
M.~Noual.
\newblock {\em Updating {A}utomata {N}etworks}.
\newblock PhD thesis, Ecole Normale Sup{\'e}rieure de Lyon-ENS LYON, 2012.

\bibitem{noual2012boolean}
M.~Noual, D.~Regnault, and S.~Sen{\'e}.
\newblock Boolean networks synchronism sensitivity and {XOR} circulant networks
  convergence time.
\newblock {\em arXiv preprint arXiv:1208.2767}, 2012.

\bibitem{noual2013non}
M.~Noual, D.~Regnault, and S.~Sen{\'e}.
\newblock About non-monotony in {B}oolean automata networks.
\newblock {\em Theoretical Computer Science}, 504:12--25, 2013.

\bibitem{park2014monomial}
J.~Park and S.~Gao.
\newblock Monomial dynamical systems in \# {P}-complete.
\newblock {\em Mathematical Journal of Interdisciplinary Sciences}, 1(1), 2012.

\bibitem{raeymaekers2002dynamics}
L.~Raeymaekers.
\newblock Dynamics of {B}oolean networks controlled by biologically meaningful
  functions.
\newblock {\em Journal of Theoretical Biology}, 218(3):331--341, 2002.

\bibitem{remy2003description}
E.~Remy, B.~Moss{\'e}, C.~Chaouiya, and D.~Thieffry.
\newblock A description of dynamical graphs associated to elementary regulatory
  circuits.
\newblock {\em Bioinformatics}, 19(suppl 2):ii172--ii178, 2003.

\bibitem{ruskey1999efficient}
F.~Ruskey and J.~Sawada.
\newblock An efficient algorithm for generating necklaces with fixed density.
\newblock {\em SIAM Journal on Computing}, 29(2):671--684, 1999.

\bibitem{ruz2013preservation}
G.~A. Ruz, M.~Montalva, and E.~Goles.
\newblock On the preservation of limit cycles in {B}oolean networks under
  different updating schemes.
\newblock {\em Advances in Artificial Life, ECAL}, pages 1085--1090, 2013.

\bibitem{sontag2008effect}
E.~Sontag, A.~Veliz-Cuba, R.~Laubenbacher, and A.~S. Jarrah.
\newblock The effect of negative feedback loops on the dynamics of {B}oolean
  networks.
\newblock {\em Biophysical Journal}, 95(2):518--526, 2008.

\bibitem{thomas1973boolean}
R.~Thomas.
\newblock Boolean formalization of genetic control circuits.
\newblock {\em Journal of Theoretical Biology}, 42(3):563--585, 1973.

\bibitem{thomas1990biological}
R.~Thomas and R.~D'Ari.
\newblock {\em Biological Feedback}.
\newblock CRC press, 1990.

\bibitem{varadarajan1990aperiodic}
K.~Varadarajan and K.~Wehrhahn.
\newblock Aperiodic rings, necklace rings, and witt vectors.
\newblock {\em Advances in Mathematics}, 81(1):1--29, 1990.

\bibitem{veliz2015dimension}
A.~Veliz-Cuba, B.~Aguilar, and R.~Laubenbacher.
\newblock Dimension reduction of large sparse and-not network models.
\newblock {\em Electronic Notes in Theoretical Computer Science}, 316:83--95,
  2015.

\bibitem{veliz2012and}
A.~Veliz-Cuba, K.~Buschur, R.~Hamershock, A.~Kniss, E.~Wolff, and
  R.~Laubenbacher.
\newblock {AND-NOT} logic framework for steady state analysis of {B}oolean
  network models.
\newblock {\em arXiv preprint arXiv:1211.5633}, 2012.

\bibitem{veliz2010dynamics}
A.~Veliz-Cuba and R.~Laubenbacher.
\newblock The dynamics of semilattice networks.
\newblock {\em arXiv preprint arXiv:1010.0359}, 2010.

\bibitem{weinblatt1972new}
H.~Weinblatt.
\newblock A new search algorithm for finding the simple cycles of a finite
  directed graph.
\newblock {\em Journal of the ACM (JACM)}, 19(1):43--56, 1972.

\end{thebibliography}



\appendix
\section{A counter example for the converse of Lemma~\ref{cycle}}    

Recall the statement of Remark~\ref{R4} that the converse of Lemma~\ref{cycle} does not hold, i.e.,  if there is a cycle of length~$m$ in each irreducible component $G_k$, it is not necessarily true that the original digraph $D$ has a cycle of length~$mp^*$. We now provide a counter example for the converse of Lemma~\ref{cycle}. 

\begin{figure}[h]
	\centering
	\includegraphics[height=60mm]{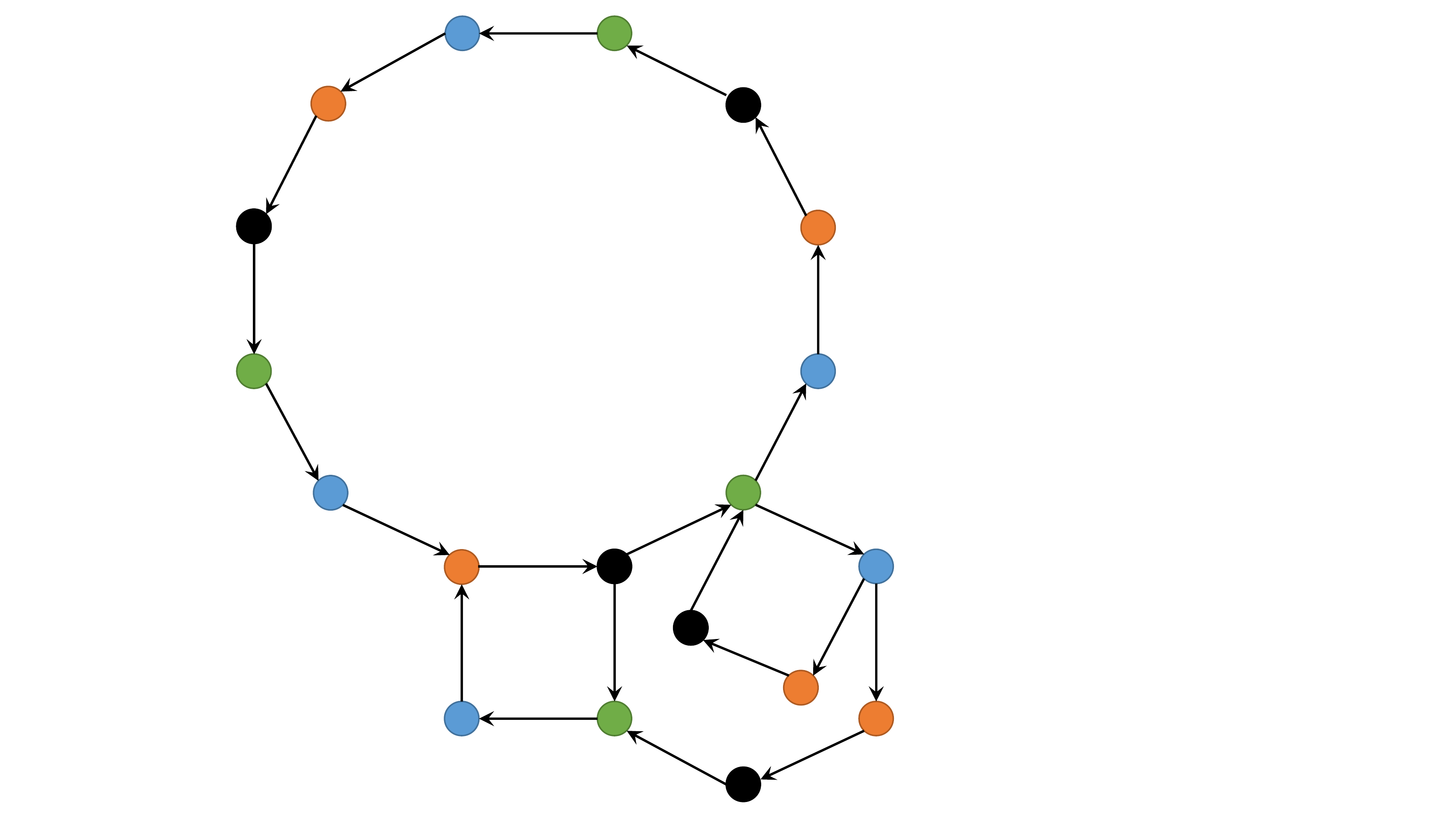}
	\caption{The digraph in the figure has four cycles whose lengths are $4$, $4$, $8$, and $12$, respectively. Let $p = 4$ be a common divisor of the cycle lengths. Then, the associated partition yields $4$ disjoint subsets, with the vertices of the same color belonging to the same subset.  
	}
	\label{counter}
\end{figure}

By slightly modifying the digraph shown in Fig.~\ref{P}: adding a new cycle of length~$4$, we obtain a new digraph shown in Fig.~\ref{counter}. This digraph has four cycles and the loop number is still $4$.

\begin{figure}[h]
	\centering
	\includegraphics[height=60mm]{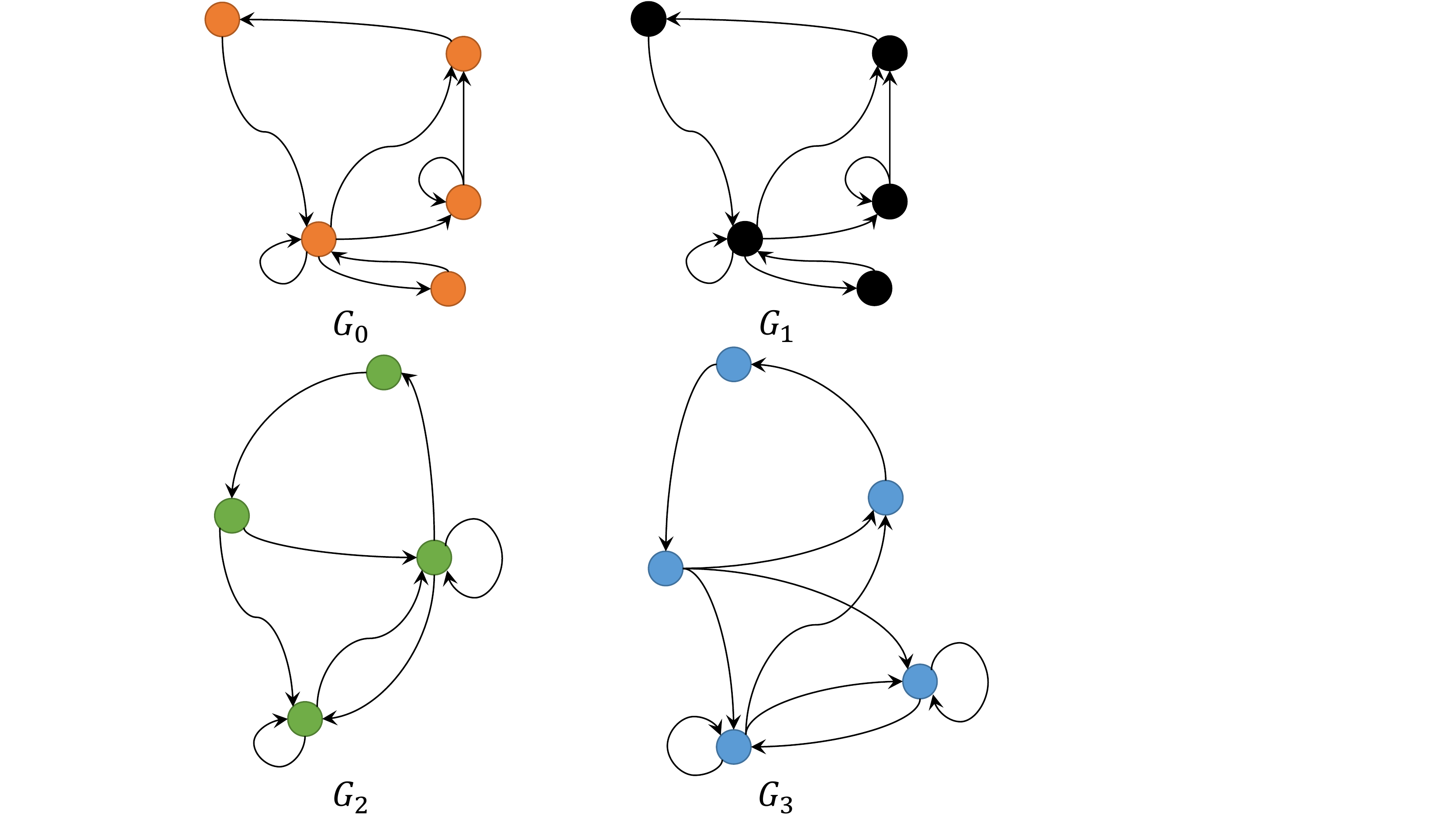}
	\caption{Irreducible components of the digraph shown in Fig.~\ref{counter} 
	}
	\label{countercomp}
\end{figure}

Following Definition~\ref{Def:IrrComp}, the irreducible components, denoted by $G_0,G_1,G_2,G_3$, are shown in Fig.~\ref{countercomp}. It can be seen that there exists a cycle of lengths~$4$ in each irreducible component. Yet, the original digraph in Fig.~\ref{counter} does not contain a cycle of length~$16$.   

\section{Proof of Lemma~\ref{equivrel}}         

We establish below the reflexivity, symmetry and transitivity of the relation ``$\sim_p$''.

{\em 1. Reflexivity.} 	
Since $D$ is strongly connected, for any $v_i\in V$, $v_i$ belongs to a cycle. Furthermore, $p$ divides the length of the cycle. We thus have $v_i\sim_p v_i$.

{\em 2. Symmetry.} 
Suppose that $v_i\sim_p v_j$; then, there exists a walk $w_{ij}$ such that $p$ divides $l(w_{ij})$. Since the graph is strongly connected, there exists a walk $w_{ji}$ from $v_j$ to $v_i$. By concatenating  $w_{ij}$ with $w_{ji}$, we obtain a closed walk $w_{ii}$ from $v_i$ to itself. It is known that any closed walk can be decomposed into cycles. This, in particular, implies that $p$ divides $l(w_{ii})$. Since $p$ divides $l(w_{ij})$, $p$ divides $l(w_{ji})$, and hence $v_j\sim_p v_i$. 

{\em 3. Transitivity.}
Suppose that $v_i\sim_p v_j$ and $v_j\sim_p v_k$; then, there exist walks $w_{ij}$ and $w_{jk}$ such that $p$ divides both $l(w_{ij})$ and $l(w_{jk})$. By concatenating $w_{ij}$ with $w_{jk}$, we obtain a walk $w_{ik}$ from $v_i$ to $v_k$. Moreover, $p$ divides $l(w_{ik})$, and hence $v_i\sim_p v_k$.
\hfill$\qed$

\section{Proof of Proposition~\ref{special}}
{\em 1. Proof of item~1.}
We first prove that if there is at least one $k=0,\ldots,p^*-1$ such that $|U_k|=1$, then $D$ is a rose. The proof is carried out by contradiction.

Suppose that the cycles in $D$ do not have the same length. Then, there exists at least one cycle whose length is greater than $p^*$. Without loss of generality, we assume that the length of $D_1$ is $rp^*$, for $r\geq 2$. Label the vertices in $D_1$ as $v_0,\ldots,v_{rp^*-1}$.  
We then define subsets $[v_0],\ldots,[v_{p^*-1}]$.   
From the second item of Proposition~\ref{417}, these subsets form a partition of $V$. Now, consider the vertices $v_{p^*}, \ldots, v_{2p^*-1}$ in $D_1$. Within $D_1$ (and hence $D$), there is a path from $v_k$ to $v_{k + p^*}$ for all $k = 0,\ldots, p^*-1$. So, we have that $v_{p^*} \in [v_0], \ldots, v_{2p^*-1} \in [v_{p^*-1}]$. In other words, $|U_k| =|[v_k]| \ge 2$ for all $k $, which contradicts the assumption that $|U_k|=1$ for at least one $k$.

We thus assume that all cycles have the same length, yet do not share a common vertex. Note that in this case, the greatest common divisor $p^*$ is the length of any cycle of $D$. Since there is at least one $k$ such that $|U_k| = 1$, we can assume, without loss of generality, that $|U_0| = 1$, and let $U_0=\{v_0\}$. Since $v_0$ is not shared by all cycles, there is a cycle $D_i$ of $D$ such that $v_i$ is not in $D_i$.  Label the vertices of $D_i$ as $v'_0,\ldots, v'_{p^*-1}$. Appealing again to the second item of Proposition~\ref{417}, we have that the subsets $[v'_0],\ldots, [v'_{p^*-1}]$ form a partition of $V$. Hence, there is some $k$ such that $v_0\in [v'_k]$. Since $v_0$ is not in $D_i$ and $v'_k$ is a vertex of $D_i$, $v_0\neq v'_k$. But then, $|U_0| =  |[v_0]| = |[v'_k]| \ge 2$, which is a contradiction.

We now prove that if $D$ is a rose, then there is at least one $k\in \{0,\ldots,p^*-1\}$ such that $|U_k|=1$. Without loss of generality, we let $v_0$ be a common vertex shared by the cycles of $D$. We now show that $|[v_0]| = 1$. Suppose not, then there is a vertex $v_i$ such that $v_i \neq v_0$ and $v_i\in [v_0]$. Since $D$ is strongly connected, $v_i$ is contained in a cycle $D_j$ of $D$. Since $v_0$ is a common vertex, $v_0$ is also contained in $D_j$. So, within the cycle $D_i$, there is a path from $v_0$ to $v_i$. Moreover, the length of the path must be greater than $0$, yet less than the length of $D_i$, which is $p^*$. On the other hand, since $v_i \sim_{p^*} v_0$, from Lemma~\ref{relation48}, the length of any walk from $v_0$ to $v_i$ has to be a multiple of $p^*$, which is a contradiction.

{\em 2. Proof of item~2.} 
If $D=(V,E)$ is a cycle digraph where the cycle has length $p^*$, then by the construction of the irreducible components, there are $p^*$ of them, each of which contains only one vertex. Now, suppose that $|U_k|=1$ for all $k=0,\ldots,p^*-1$; we show that $D$ is a cycle digraph where the cycle has length $p^*$. Let vertices $v_0,\ldots,v_{p^*-1}$  be such that $U_k = \{v_k\}$ for all $k = 0,\ldots, p^*-1$. Since $U_0,\ldots, U_{p^*-1}$ form a partition of $V$, we have that the vertex set of $D$ is given by $\{v_0,\ldots, v_{p^*-1}\}$. Since the number of vertices of $D$ is $p^*$, the length of any cycle of $D$ is no greater than $p^*$. Furthermore, if the loop number of $D$ is $p^*$, then each cycle of $D$ has to be a Hamiltonian cycle (i.e., a cycle that passes all the vertices of $D$). But this happens if and only if $D$ itself is a cycle digraph.  
\qed
                                        
\end{document}